\documentclass[sn-mathphys]{sn-jnl}

\jyear{2021}%

\usepackage{amssymb}
\usepackage{amsthm}
\usepackage{amsmath}
\usepackage{amsfonts}
\usepackage{hyperref}
\usepackage{color}

\usepackage[T1]{fontenc}
\usepackage[latin2]{inputenc}

\theoremstyle{thmstyleone}%

\newtheorem{theo}{\bf Theorem}[section]

\newtheorem{lem}{\bf Lemma}[section]
\newtheorem{rem}{\bf Remark}[section]
\newtheorem{defi}{\bf Definition}[section]

\newtheorem{fact}{\bf Fact}[section]

\newtheorem{oq}{\bf Open Question}[section]

\newcommand{\rp}{\mathbb{R}_+}

\theoremstyle{thmstyletwo}%

\theoremstyle{thmstylethree}%

\newcommand{\barint}{
         \rule[.036in]{.12in}{.009in}\kern-.16in
          \displaystyle\int  }

\begin{document}

\title[Asymptotics, trace, and density results for  weighted Dirichlet  spaces]{Asymptotics, trace, and density results for  weighted Dirichlet  spaces defined on the halfline}


\author[1]{\fnm{} \sur{Claudia Capone}}\email{c.capone@na.iac.cnr.it}
\equalcont{These authors contributed equally to this work.}

\author*[2]{\fnm{} \sur{Agnieszka Ka\l{}amajska}}\email{A.Kalamajska@mimuw.edu.pl}
\equalcont{These authors contributed equally to this work.}


\affil*[1]{\orgdiv{Istituto per le Applicazioni del Calcolo ``Mauro Picone''}, \orgname{Consiglio Nazionale delle Ricerche}, \orgaddress{\street{via Pietro Castellino 11}, \city{Napoli}, \postcode{80131}, 
\country{Italy}}}

\affil[2]{\orgdiv{Faculty of Mathematics, Informatics, and Mechanics}, \orgname{University of Warsaw}, \orgaddress{\street{ul.
Banacha 2}, \city{Warsaw}, \postcode{02--097}, 
\country{Poland}}}



\abstract{We give analytic description for the completion of $C_0^\infty (\rp)$ in Dirichlet space
$D^{1,p}(\rp, \omega):=
\{  u:\rp\rightarrow \mathbf{R}:  u\ \hbox{\rm is locally absolutely }  $\\
$ \hbox{\rm  continuous on}\ \rp \ \hbox{\rm and} \  \| u^{'}\|_{L^p(\rp, \omega)}<\infty \}$, for  given continuous positive weight $\omega$ defined on $\rp$, where $1<p<\infty$. The conditions are described in terms of the modified variants of the  $B_p$ conditions due to Kufner and Opic from 1984, which in our approach are focusing on integrability of $\omega^{-p/(p-1)}$ near zero or near infinity. 
Moreover, we propose applications of our results to: obtaining new variants of Hardy inequality,  interpretation of boundary value problems  in ODE's defined on the halpfline with solutions in $D^{1,p}(\rp, \omega)$, new results from complex interpolation theory dealing with interpolation spaces between weighted Dirichlet spaces,  and  to
 derivation of new Morrey type embedding theorems for our Dirichlet space.}





\maketitle

\newpage
\section{Introduction}\label{sec1}
In this paper we are interested in  weighted Dirichlet spaces
\begin{eqnarray*}\label{przestrzen}
D^{1,p}(\rp, \omega)=\{  u:\rp\rightarrow \mathbf{R}:  u\  \hbox{\rm is locally absolutely continuous on}\ \rp \ {\rm and} \\ \| u^{'}\|_{L^p(\rp, \omega)}<\infty \}.\nonumber
\end{eqnarray*}
In most situations we assume that the weight {$\omega :\rp \to \rp$}, is continuous  and
$1<p<\infty$.

In some cases we also assume that $\omega$ satisfies
the  localized at the endpoint variant  of the general $B_p$-condition due to Kufner and Opic from \cite{kuf-opic}:

$$
 B_p(0):\ \int_{(0,1)}\omega(t)^{-1/(p-1)}dt <\infty \ \ {\rm or}\  \
 B_p(\infty):\ \int_{(1,\infty)}\omega(t)^{-1/(p-1)}dt <\infty .
$$

\noindent We address and analyze several problems related to such spaces.

\smallskip
{\it Asymptotic behaviour.} One of the topics of our interest is  asymptotic behaviour at the endpoints  for elements of such spaces.
Assume for example that $\omega\in B_p(0)$.
Among our results in this direction, we show in Theorem \ref{r-sets} that
when $\omega\in B_p(0)$, $u\in D^p(\rp,\omega)$ and $c\in \mathbf{R}$,  then the conditions (a),(b),(c) are equivalent, where
\begin{eqnarray}
(a) && \exists_{t_n\, \searrow \, 0}:\lim_{t_n\to 0}(u(t_n)-c)
= 0,\ \ \
(b) \ \ \
\lim_{t\to 0} u(t) =c, \label{granice}\\
(c) &&
 \lim_{t\to 0}\frac{u(t)-c}{\Omega^0_\omega (t)}= 0 \,, \, \ {\rm where}\ \, \displaystyle{\Omega_{\omega}^0 (t):= \left( \int_{0}^t \omega (\tau )^{\frac{-1}{p-1}  }d\tau\right)^{1-\frac{1}{p}}}.\nonumber
\end{eqnarray}
As $\Omega_\omega^0 (t)\to 0$ when $t\to 0$, we
clearly have $(c)\Rightarrow (b) \Rightarrow (a)$. The nontrivial part {is} to prove that the converse implications hold.

\noindent Similar analysis is also provided about behaviour near infinity.

\smallskip
{\it Trace operator.}
There are several ways to define the trace of Sobolev function, see e.g. \cite{kjf}, Section 6.10.5 for the classical approach.
We ask about the limit (b) in (\ref{granice}) and  we define
\begin{equation}\label{traceop}
Tr^{0}u:=\lim_{t\to 0} u(t).
\end{equation}
Clearly, one has to ask if such limit is well prescribed in our Dirichlet space setting. It is always so, when we assume
that $\omega\in B_p(0)$, see Theorem \ref{r-sets}, part iii) and it is never so, when $\omega\not\in B_p(0)$, see Theorem \ref{domknietosczbiorowzero}.

\smallskip
{\it The norm on Dirichlet space.} { Let us note that the quantity $\| u^{'}\|_{L^p(\rp, \omega)}$ annihilates all constants,} therefore $\| u^{'}\|_{L^p(\rp, \omega)}$ cannot define the norm on Dirichlet space $D^{1,p}(\rp,\omega)$.
However,  for any $a\in\rp$, the quantity
\begin{equation}
\| u\|_{D^{1,p}(\rp,\omega)}^{(a)}:=  \lvert u(a) \rvert + \| u^{'}\|_{L^p(\rp, \omega)}, 
\label{normaa}
\end{equation}
defines the  norm on $D^{1,p}(\rp,\omega)$ and makes it a Banach space. Moreover,  all such norms $\|\cdot\|_{D^{1,p}(\rp,\omega)}^{(a)} : a\in\rp$,  are equivalent. See Fact \ref{norm-dirichlet}.

In the case of $\omega\in B_p(0)$, we can extend the definition of the norm \eqref{normaa} also
to  $a=0$, by  putting
$u(0):= Tr^{0}(u)$ in place $u(a)$. Such modification gives also the equivalent norm. In such case the trace operator is continuous as functional on our Dirichlet space  equipped with any of the proposed norms (\ref{normaa}), including $a=0$, see Theorem \ref{r-sets}, part iv).
In case of $\omega\in B_p(\infty)$, similar property holds with
$a=\infty$, {  see Theorem \ref{r-setsinfty}}.

\smallskip
{\it Representation of functions.} Let us focus on the case  of $\omega\in B_p(0)$.
Because in that case  the limit $\displaystyle{\lim_{t\to 0}u(t)}$ does exist  for any $u$ in Dirichlet space and $u^{'}$ is integrable near zero, {every element   $u\in D^{1,p}(\rp,\omega)$}
can be represented as
$$
u(t)-Tr^{0}u = \int_0^t u^{'}(\tau )d\tau =Hu^{'}(t),
$$
 where {{on the}} right hand side {above} we deal with Hardy transform of $u^{'}$,
remembering that $u^{'}$  belongs to $L^p(\rp,\omega)$. This allows to deduce several further properties,  for example applications to Hardy inequality, see Section \ref{hardyin}.
{Similar representations hold in case  of $\omega\in B_p(\infty)$.
In this case we use the conjugate Hardy transform as in (\ref{hardypresentation}), see Theorem \ref{r-setsinfty}.
}

\smallskip
{\it Questions about densities.} Let us denote  by $D_0^{1,p}(\rp, \omega)$ the completion of $C_0^{\infty}(\rp)$  in $D^{1,p}(\rp, \omega)$ in any norm like (\ref{normaa}).
It is the  natural question to ask about characterization of weights $\omega$, for which
$D_0^{1,p}(\rp, \omega) = D^{1,p}(\rp, \omega)$. If that spaces are not the same, we can ask if it is possible to characterize completely the space $D_0^{1,p}(\rp, \omega)$ by some analytic conditions expressed in
terms of the weight $\omega$. Let us focus {again} on the case of
$\omega\in B_p(0)$.  {In Theorem \ref{characterization} we have proved
that for such $\omega$}
\begin{eqnarray*}
 D_0^{1,p}(\rp,\omega )= {\mathcal B}^{0}_{p,\omega}(0):=\{ u\in D^{1,p}(\rp,\omega):  \lim_{t\to 0}u(t)=0\}
 \Longleftrightarrow \omega\not\in B_p(\infty).
\end{eqnarray*}
This gives the analytic  characterization of weights for which $ D_0^{1,p}(\rp,\omega )$ is the kernel of trace operator as in \eqref{traceop}, in the case of $\omega\in B_p(0)$. Let us emphasize that the trace operator  $u\mapsto Tr^{0}u$  in such case is well defined and continuous.

{As the consequence of Theorem \ref{r-sets} and Theorem \ref{domknietosczbiorowzero}, in the case of $\omega\not\in B_p(0)$, the trace operator $Tr^0(\cdot)$ at zero is not well defined.}
We can thus ask question if the space $D_0^{1,p}(\rp,\omega )$  could {still} be characterized by some analytic conditions, without assuming that $\omega\in B_p(0)$. Such characterization is provided in Theorem \ref{glowne}, which gives the precise analytic characterization of $D_0^{1,p}(\rp,\omega )$, expressed  in terms of the conditions $B_p(0)$ and $B_p(\infty)$. For example, as follows from Theorems: \ref{glowne} and \ref{zamkniecieimplikacji}, among the other statements, we show that
\begin{eqnarray*}
 D_0^{1,p}(\rp,\omega ) = {\mathcal B}^{0}_{p,\omega}(0)
\Longleftrightarrow \omega \in B_p(0)\setminus  B_p(\infty).
\end{eqnarray*}

\bigskip
{\it Applications to: Hardy inequality, Boundary Value Problems
 (B.V.P.) in ODE's, generalized Morrey  Theorem and to complex interpolation theory.}
Having more precise information about representation of function from  our Dirichlet space, or about the  asymptotic behaviour of the functions from given Dirichlet space near zero or infinity, one can deduce
more precise variants of Hardy and  conjugate Hardy inequality (see Section \ref{hardyin}), or establish
 if the given boundary value problem, presented in term of vanishing of  function near zero or infinity in the analyzed ODE, is well posed or not.  The discussion is provided in Section \ref{bvpb}. {Moreover, in Section \ref{morrey-general} }  we focus  on certain generalization of  Morrey  Theorem, which deals with $B_p$-conditions. In Section \ref{interpo} we have also presented some new applications of our results to complex  interpolation
 theory,  dealing with weighted Dirichlet type spaces,  inspired by questions posed recently in \cite{cwi-ein}.

\bigskip
{\it Novelty and link with literature.}
To our best knowledge, our results concerning the {\it asymptotic behaviour} near the endpoints of the interval of functions in the non trivially weighted Dirichlet spaces, as summarized in Theorems \ref{r-sets} and \ref{r-setsinfty}, are new. In the non-weighted setting they are motivated by Morrey theorem, see Section \ref{morrey-general}. However, similar type conclusions can be found also in the case of power weight in \cite{akkppstudia09}, on page 9.

{\it Density} results in the general Dirichlet space setting,
are rather missing in the literature.
In the case
of $\omega\equiv 1$, they were obtained first by Sobolev in 1963 (\cite{sobolev}) and now they are well understood. See also e.g. \cite{haj-kalam}, Theorem 4 and references therein, where density results are obtained with respect to the semi norm $\| u^{'}\|_{L^p(\rp)}$  as in Fact \ref{norm-dirichlet},
instead of the norm $\| u\|^{(a)}$ from \eqref{normaua}.
Our density anaylisis is based on the localized at endpoints $B_p$ conditions as in Definition  \ref{bpcondi},
which were not considered before. However,
some preliminary ideas for such conditions can be found in \cite{kuf-opic}, see Remark \ref{kuopbp}.

Most of the classical density results  deal with Sobolev spaces,  not Dirichlet spaces.
In case of Sobolev spaces the additional restriction on function $u$ is provided, that is its integrability with some power.
We would like to emphasise that our density results mostly deal with the norm (\ref{normaa}), they are restricted to Dirichlet (not Sobolev) space, and  characterize completely the admitted  weights.

As about {\it the tools}, for the analysis of asymptotic behaviour
we use  simple computations based on Taylor's formula in $1$-d.
To study density, we propose the technique, which in our opinion is new in such setting.
 We call it the energy - caloric approximation, as it is based on the variational technique. More precisely, we first find the function
which on segments $[a,b]$ minimize the energy functional $$\int_{(a,b)}
{  \lvert u^{'}(\tau )  \rvert }^{p} \omega (\tau)d\tau   $$
 with the given boundary data at the endpoints $\{ a, b\}$. In further step we
extend such local minimizers to compactly supported functions in the same Dirichlet-Sobolev class.  See the considerations in Sections \ref{sec-density-0} and \ref{sec-density-infty}.

{\it Organization of the paper.}  After the preliminary
results presented in Section \ref{preliminaries}, we analyze questions about the asymptotic behaviour and trace in Section \ref{asympto}, while
density results are presented in Section \ref{densityres}. Main applications: to the derivation of Hardy inequality, to the well posedness of B.V.P., to the derivation of Morrey type theorems, as well as to complex interpolation theory in Dirichlet space setting,  are discussed in
Section \ref{appli}. Some additional remarks
are presented in Section \ref{furtherdev}, while
in Section \ref{appen} we enclose some auxiliary computations and complementary results, for reader's convenience.

\section{Notation and Preliminaries}\label{preliminaries}

\subsection{Basic notation}

\noindent In most situations we deal with  positive continuous functions
$\omega:\rp\rightarrow \rp$, referred as positive weights, where, by positive expression, we mean that it is strictly larger than zero.
However for our purposes we  consider continuous  weights only,
we will sometimes formulate our statements in the more general setting.

We use standard notation: $C_0^\infty(\rp)$, $L^p(\rp,\omega)$,
$L^p_{loc}(\rp,\omega)$, $Lip (\rp)$, $W^{1,p}(\rp)$, $W^{1,p}_{loc}(\rp)$, for smooth compactly supported functions, weighted $L^p$-spaces and their local variants, Lipschitz functions, the classical Sobolev spaces and their local variants. {We} will also use the more specific notation for the local variants of $L^p$ and Sobolev - type spaces.
For $1\le p<\infty$, by $L^p_{loc}([0,\infty))$  we denote all functions
$f\in  L^p_{loc}(\rp)$ which are $p$-integrable near zero
(shortly
$\int_0 {\lvert f \rvert}^p d\tau <\infty$), while by $L^p_{loc}((0,\infty])$  we denote all functions
$f\in  L^p_{loc}(\rp)$, which are $p$-integrable near infinity
(shortly
$\int^\infty { \lvert f \rvert}^p d\tau <\infty$).
Analogous definitions with obvious modifications will be used to denote the corresponding Sobolev spaces: $W^{1,p}_{loc}([0,\infty))$, $W^{1,p}_{loc}((0,\infty ])$, and their generalizations.
 In most situations we will refer to the Lebesgue integral. However,   sometimes we will also refer to the Newtonian interpretation of the integral when writing $\int_a^bf dx (=F(b)-F(a))$ where $F^{'}=f$ a.e., in place of $\int_{(a,b)}f dx$. By  measurable sets, we mean sets that are measurable with respect to the Lebesgue measure.

\noindent

{Let $X$ be some subset of Lebesgue measurable functions defined on $\rp$. By $X_c$ we will  denote its subset  consisiting of functions with compact suport in $\rp$.
When $X\subseteq Z$, where $(Z,\|\cdot\|_Z)$ is  some Banach space}, then by $\overline{X}^{\|\cdot\|_Z}$ will denote the completion of $X$ in  the norm  $\| \cdot\|_Z$.
The symbol  $Z_0$ will be reserved for $\overline{(C_0^\infty(\rp)\cap Z)}^{\|\cdot\|_Z}$.

In our estimates, we will sometimes write $f\sim 1$ if the  function $f$ defined on its respective domain  can be estimated from both sides  by positive {constants}, while the notation $f\precsim 1$ will mean that the function is bounded from above.

\bigskip
\noindent
\subsection{General and local $B_p$-conditions for weights}
We will deal with the following variants of the $B_p$-condition introduced by Kufner and Opic in \cite{kuf-opic}.

\begin{defi}[$B_p$-conditions]\label{bpcondi}\rm  Let $\omega: \rp \rightarrow {[0,\infty)}$ be a measurable function which is positive almost everywhere, $1<p<\infty$. We say that
\begin{description}
\item[a)]
   $\omega$ is a $B_p$-weight (shortly $\omega\in B_p$)  if
$
\omega^{-1/(p-1)}\in L^1_{{\rm loc}}((0,\infty)),
$ see \cite{kuf-opic};
\item[b)]
  $\omega$ is a $B_p$-weight near zero (shortly $\omega\in B_p(0)$)  if $
\displaystyle{\int_0\omega^{-1/(p-1)}d\tau <\infty;}
$
\item[c)]
  $\omega$ is a $B_p$-weight near infinity (shortly $\omega\in B_p(\infty)$)  if
   $\displaystyle{\int^\infty\omega^{-1/(p-1)}d\tau <\infty .}$
\end{description}
\end{defi}

\noindent Note that both conditions $B_p(0)$ and $B_p(\infty)$ imply that $\omega \in B_p$. Moreover, by H\"older inequality, for any measurable set $K\subseteq [0,\infty)$
\begin{eqnarray}
\label{bp-rachunek}
\int_{K}  \lvert
f (\tau)  \rvert d\tau =\nonumber\\
\int_K  \lvert f(\tau)  \rvert\omega(\tau)^{\frac{1}{p}}
\omega(\tau)^{-\frac{1}{p}} \le \left( \int_K  {\lvert f(\tau) \rvert}^p\omega(\tau)d\tau \right)^{\frac{1}{p}} \left( \int_K \omega(\tau)^{\frac{-1}{p-1}}d\tau \right)^{1-\frac{1}{p}},
\end{eqnarray}
This implies.

\begin{fact}\label{zanurzenia}
Let $\omega, p$ be as in Definition \ref{bpcondi}.
The following statements hold:
\begin{description}
\item[i)] when $\omega\in B_p$ then $L^p(\rp,\omega)\subseteq L^1_{loc}(\rp)$;
\item[ii)] when $\omega\in B_p(0)$ then $L^p(\rp,\omega)\subseteq L^1_{loc}([0,\infty))$;
\item[iii)]
when $\omega\in B_p(\infty)$ then $L^p(\rp,\omega)\subseteq L^1_{loc}((0,\infty])$.
\end{description}
\end{fact}

In our specific situation, we assume that the weight $\omega$ is continuous and positive, which guarantees that $\omega\in B_p$.
{The conditions} $B_p(0)$ and $B_p(\infty)$, which to our best knowledge were not introduced earlier, are motivated by the general $B_p$ condition from \cite{kuf-opic}. More precise information about $B_p$-conditions is provided in Remark \ref{kuopbp}.

\smallskip

\subsection{Weighted Dirichlet spaces}

\noindent
We start with the definition of weighted Dirichlet space.

\begin{defi}[weighted Dirichlet space]\label{trutututu}\rm
Let $\omega : \rp \rightarrow [0,\infty)$ be positive weight, that is $\omega >0$\  a.e., $1<p<\infty$.
By $D^{1,p}(\rp ,\omega)$  we will denote the Dirichlet space consisting with all functions $u\in W^{1,1}_{loc}(\rp)$ such that
\begin{eqnarray*}
\| u\|^{*}_{D^{1,p}(\rp ,\omega)}:= \left(
\int_{\rp} { \lvert u^{'}(t) \rvert }^p\omega (t)dt\right)^{\frac{1}{p}} <\infty  .
\label{weight-dirichlet}
\end{eqnarray*}
\end{defi}
Clearly, the expression $\| u\|^{{*}}_{D^{1,p}(\rp ,\omega)}$  annihilates constant functions, so it defines the semi norm on $D^{1,p}(\rp ,\omega)$ but not the norm.

\smallskip
\noindent
We are interested in Dirichlet spaces in the case when
$\omega$ is continuous and positive, and so
$\omega\in B_p$. In that case we show
that the homogeneous Dirichlet space $\tilde{D}^{1,p}(\rp,\omega)$ defined below  is complete.
The proof is enclosed in the Appendix for reader's convenience.

\begin{fact}[homogeneous Dirichlet space]\label{norm-dirichlet}
Let $\omega, p$ be as in Definition \ref{trutututu}  and  consider the relation:
$u\sim v$ when $u,v\in D^{1,p}(\rp ,\omega)$ and $u-v\equiv c$ for some constant $c\in\mathbf{R}$. Then
\begin{eqnarray*}
\tilde{D}^{1,p}(\rp ,\omega)&:=& {D^{1,p}(\rp ,\omega)}\big{/}{\sim } \ \ \hbox{\rm equipped with the norm}\\
 \| \{ u+c\}_{c\in\mathbf{R}}\|^{*}_{ {\tilde{D}^{1,p}(\rp ,\omega)}}&:=&
\| u^{'}\|_{L^p(\rp,\omega)}
\end{eqnarray*}
 is a Banach space.
\end{fact}

\noindent In the following fact we analyze the norms in Dirichlet spaces.

\begin{fact}[the norms on $D^{1,p}(\rp ,\omega)$ and topology of convergence]\label{bp-dirichlet}
 Let $\omega, p$ be as in Definition \ref{trutututu}. Then for any $a\in (0,\infty)$ the
expression
\begin{equation}\label{normaua}
\| u\|_{D^{1,p}(\rp,\omega)}^{(a)}:= \| u^{'}\|_{L^p(\rp,\omega)} +  \lvert  u(a) \rvert 
\end{equation}
is the norm on $D^{1,p}(\rp ,\omega)$, which makes $D^{1,p}(\rp ,\omega)$ a Banach space.

Moreover, for all $a\in \rp$ the norms $\|\cdot\|_{D^{1,p}(\rp,\omega)}^{(a)}$ are equivalent on $D^{1,p}(\rp ,\omega)$ and
\begin{eqnarray}\label{charaktetyzacjazb}
\|u_n-u\|_{D^{1,p}(\rp,\omega)}^{(a)}\stackrel{n\to\infty}{\rightarrow} 0 \Longleftrightarrow ~~~~~~~~~~~~~~~~~~~~~~~~~~~\\
 \left( u_n^{'}\to u^{'}\
\hbox{\rm  in}\  L^p(\rp,\omega)\ \hbox{\rm  and} \
u_n\to u\ \hbox{\rm uniformly on compact sets in}\  \rp \right) .\nonumber
\end{eqnarray}
\end{fact}

\smallskip
\noindent
{\bf Proof.} We observe that $D^{1,p}(\rp ,\omega)\subseteq W^{1,1}_{loc}(\rp)\subseteq C(\rp)$ and so the value
$u(a)$ is well prescribed. In particular $\|\cdot\|_{D^{1,p}(\rp,\omega)}^{(a)}$ is the norm on $D^{1,p}(\rp ,\omega)$.

\noindent {Moreover, $(D^{1,p}(\rp,\omega),\| \cdot\|^{(a)}_{D^{1,p}(\rp,\omega)})$}
is a Banach space, because when $\{ u_n\}$ is the Cauchy sequence in $D^{1,p}(\rp ,\omega)$, then, due to Fact
\ref{norm-dirichlet}, there exists $v\in D^{1,p}(\rp ,\omega)$ such that
$u_n^{'}$ converge to $v^{'}$ in $L^p(\rp,\omega)$.

\noindent Then
for $u(t):=\displaystyle{\int_a^t v^{'}(\tau)d\tau + \lim_{n\to\infty} u_n(a)}$ we have
$\| u_n-u\|_{D^{1,p}(\rp,\omega)}^{(a)}\to 0$ as $n\to\infty$.
The equivalence of norms is a consequence of the following estimate holding for any closed interval $I$ such that $a,b\in I\subset \rp$, $0<a<b<\infty$:
\begin{eqnarray}
 \lvert u(b)- u(a)  \rvert\le \int_{(a,b)}  \lvert u^{'}(\tau)  \rvert\omega (\tau)^{\frac{1}{p}} \omega (\tau)^{-\frac{1}{p}} d\tau \nonumber\\
\le
\left( \int_I  {\lvert u^{'}(\tau)  \rvert}^p\omega d\tau\right)^{\frac{1}{p}} \left( \int_I \omega (\tau)^{\frac{-1}{p-1}}d\tau\right)^{1-\frac{1}{p}}
= C_I \left( \int_0^\infty { \lvert u^{'}(\tau)  \rvert}^p\omega d\tau\right)^{\frac{1}{p}}
,\label{czwartek}
\end{eqnarray}
where $C_I:= \displaystyle{\left( \int_I \omega (\tau)^{\frac{-1}{p-1}}d\tau\right)^{1-\frac{1}{p}}}$.
As a consequence of  (\ref{czwartek}) we get (\ref{charaktetyzacjazb}) with any $b\in I$ and  $u_n-u$ in place of $u$.\hfill$\Box$

\bigskip

More precise analysis, dealing with the conditions $B_p(0)$ and $B_p(\infty)$, will be provided in our next section.

\section{{Asymptotics and trace}
}\label{asympto}

\subsection{Analysis in the case of $\omega \in B_p(0)$}

We start with the analysis within the case of $\omega\in B_p(0)$.
We obtain the following statement, which deals with  trace operator defined  at zero {$Tr^{0}(\cdot)$}, as in \eqref{trzero} below,
 and precisely describes the elements of weighted Dirichlet space.

\begin{theo}[asymptotic behaviour and trace at zero]\label{r-sets}
Let $\omega :\mathbf{R}_+\rightarrow \mathbf{R}_+$,\\  $\omega\in C(\rp)\cap  B_p(0)$, $1<p<\infty$ and consider the following subsets in
$D^{1,p}(\rp , \omega)$, defined for $c\in \mathbf{R}$:
\begin{eqnarray}\label{sladyzero}
{\mathcal R}^{0}_{p,\omega}(c)&:=&\left\{ u\in D^{1,p}(\rp , \omega),  u=\int_0^t u^{'}(\tau)d\tau +c : u^{'}\in L^p(\rp,\omega)\right\}; \\
{\mathcal A}_{p,\omega }^{0}(c)&:=&\left\{ u\in D^{1,p}(\rp , \omega) :   \exists_{t_n \, \searrow \, 0}: \lim_{n\to \infty} u(t_n)=c \right\};\nonumber\\
{\mathcal B}_{p,\omega }^{0}(c)&:=&  \left\{ u\in D^{1,p}(\rp , \omega ) : \lim_{t\to 0} u(t) =c\right\} ;\nonumber\\
{\mathcal C}^{0}_{p, \omega}(c)  &:=&  \left\{ u\in D^{1,p}(\rp , \omega) :
 \lim_{t\to 0} \frac{u(t)-c}{\Omega_{\omega}^0 (t)} =0,\ \  {\rm sup_{t>0}}\, \frac{u(t)-c}{\Omega_{\omega}^0 (t)}<\infty ,\right\}, \ \nonumber
\end{eqnarray}

\noindent where $\displaystyle{{\Omega}_{\omega}^0 (t)}= \left(
 \int_{0}^t \omega (\tau )^{\frac{-1}{p-1}} d\tau
 \right)^{1-\frac{1}{p}}$.

\medskip
\noindent {The following statements hold.}
\begin{description}
\item[i)] For any $c\in \mathbf{R}$, the set ${\mathcal R}^{0}_{p,\omega}(c)$ is a  closed subset in $D^{1,p}(\rp,\omega)$, equipped with any norm $\|\cdot\|^{(a)}$ as in
\eqref{normaua}, where $a\in\rp$.
\item[ii)]  For any $c\in \mathbf{R}$
\begin{eqnarray}\label{rdef}
{\mathcal A}_{p,\omega }^{0}(c)={\mathcal B}_{p,\omega }^{0}(c)={\mathcal C}^{0}_{p, \omega}(c)= {\mathcal R}^{0}_{p,\omega}(c).
\end{eqnarray}
\item[iii)]
 For every $u\in D^{1,p}(\rp,\omega)$ there is $c\in \mathbf{R}$ such that $u\in {\mathcal R}^{0}_{p,\omega}(c)$. In particular, the trace operator
\begin{equation}\label{trzero}
Tr^{0}(u):= \lim_{t\to 0} u(t)=: u(0),
\end{equation}
is well defined for every $u\in D^{1,p}(\rp, \omega)$ and $D^{1,p}(\rp ,\omega) = \cup_{c\in\mathbf{R}}  {\mathcal R}^{0}_{p,\omega}(c)$.

\noindent {Moreover,} 
\begin{equation*}\label{dirrepp}
D^{1,p}(\rp,\omega )  = \left\{ u= \int_0^t v(\tau)d\tau + c  : c\in\mathbf{R},  v\in L^p(\rp,\omega)\right\}.
\end{equation*}
\item[iv)] The quantity
\begin{equation}\label{normazero}
 \| u\|_{D^{1,p}(\rp,\omega)}^{(0)}: =
 \| u^{'}\|_{L^p(\rp,\omega)}+  \lvert  Tr^{0}(u) \rvert 
 \end{equation}
is the norm on $D^{1,p}(\rp,\omega )$, which is
equivalent to any norm $\| u\|^{(a)}_{{D^{1,p}(\rp,\omega)}}$, where
$a\in \rp$.
\end{description}
\end{theo}

\noindent {\bf Proof.}
We  observe that the substitution of $u-c$ in place of $u$,  reduces the proofs of i) and ii)
 to the case of $c=0$. Therefore, for that statements, we only show the case of $c=0$.

\smallskip
\noindent
{\bf i):} \, Let us consider a sequence $\{u_n\}_n \subseteq {\mathcal R}^{0}_{p,\omega} (0)$, $u_n(t)\stackrel{n\to\infty}{\to} u(t)$ in $D^{1,p}(\rp,\omega)$. {By} Fact \ref{bp-dirichlet}
$$u_n \to u  \quad {\rm uniformly \, on \, compact \, sets\,  in \,} \rp\ \
{\rm and}\ \
u_n^{'} \to u^{'} \ {in}\   L^p(\rp,\omega) \, .$$

\noindent As $\omega \in B_p(0)$, therefore the above convergence yields $u_n^{'} \to u^{'}$ in $L^1(0,K)$, for every $K>0$.
Since $u_n \in {\mathcal R}^{0}_{p,\omega}(0)$, for every $t \, {>0}$, we have
$$u_n(t)=\int_0^t u_n^{'}(\tau) d\tau \stackrel{n\to\infty}{\rightarrow} u(t)=\int_0^t u^{'}(\tau) d\tau .$$
In particular  $u\in {\mathcal R}^{0}_{p,\omega}(0)$ and so ${\mathcal R}^{0}_{p,\omega}(0)$ is closed.
\smallskip

\noindent {\bf ii):} \, We start by proving the identity $ {\mathcal R}^{0}_{p,\omega}(0)= {\mathcal B}_{p,\omega }^{0}(0) $.

 \noindent Let $u\in  {\mathcal R}^{0}_{p,\omega}(0)$. Then
\begin{equation}\label{ajeden}
u(t)=\int_0^t u^{'}(\tau) d\tau .
\end{equation}
Therefore  $u(t) \to 0$ as $t \to 0$, hence $u \in {\mathcal B}_{p,\omega }^{0}(0)$. {This} gives ${\mathcal R}^{0}_{p,\omega}(0)\subseteq {\mathcal B}_{p,\omega }^{0}(0)$.

\noindent On the other hand, if $u \in {\mathcal B}_{p,\omega }^{0}(0)$, then for every $0<\overline{t}<t<K$, we have
$$u(t)-u(\overline{t})= \int_{\overline{t}}^t u^{'}(\tau) d\tau .$$
Since $u^{'} \in L^1((0,K))$ for any  $K$, by taking the limit as $\overline{t} \to 0$, we get \eqref{ajeden}.
Hence $u \in {\mathcal R}^{0}_{p,\omega}(0)$.

\smallskip
\noindent
 We will complete the proof of ii) by proving the identity: ${\mathcal A}_{p,\omega }^{0}(0)={\mathcal B}_{p,\omega }^{0}(0)={\mathcal C}_{p,\omega }^{0}(0)$.

\noindent We first show the equality
\begin{equation}\label{a}
{\mathcal A}_{p,\omega }^{0}(0)={\mathcal B}_{p,\omega }^{0}(0).
\end{equation}
Clearly, ${\mathcal A}_{p,\omega }^{0}(0)\supseteq {\mathcal B}_{p,\omega }^{0}(0)$, so we have to prove the converse inclusion.
To this aim, let us take $u\in {\mathcal A}_{p,\omega }^{0}(0)$, and let $t_n\to 0$ be such that $u(t_n)\to 0$ as $n\to \infty$. Then,  for any $t$ such that $0<t_n<t$:
\begin{eqnarray}
 \lvert u(t)-u(t_n)  \rvert &\le& \int_{t_n}^t  \lvert u^{'}(\tau )  \rvert d\tau = \int_{t_n}^t  \lvert u^{'}(\tau )  \rvert \omega^{\frac{1}{p}}(\tau )\omega^{-\frac{1}{p}}(\tau ) d\tau \label{est1}\\
&\le&
\left( \int_{t_n}^t  {\lvert u^{'}(\tau )  \rvert}^p\omega (\tau )d\tau\right)^{\frac{1}{p}} \left( \int_{t_n}^t \omega (\tau )^{\frac{-1}{p-1}  }d\tau\right)^{1-\frac{1}{p}} \nonumber\\
&\le&
\left( \int_{0}^t   {\lvert u^{'}(\tau )  \rvert}^p\omega (\tau )d\tau\right)^{\frac{1}{p}} \left( \int_{0}^t \omega (\tau )^{\frac{-1}{p-1}  }d\tau\right)^{1-\frac{1}{p}} . \nonumber
\end{eqnarray}

\noindent After letting  $n\to\infty$, we get
\begin{eqnarray}\label{b}
 \lvert u(t)  \rvert\le \left( \int_{0}^t  {\lvert u^{'}(\tau )  \rvert}^p\omega (\tau )d\tau\right)^{\frac{1}{p}} \cdot \Omega_{\omega}^0(t) 
\stackrel{t\to 0}{\rightarrow}0,
\end{eqnarray}
which proves (\ref{a}).

\smallskip
\noindent
Let us show that both sets in (\ref{a}) are the same as ${\mathcal C}^{0}_{p, \omega}(0)$. Indeed, let us consider $u\in {\mathcal A}_{p,\omega }^{0}(0)$. Then, by (\ref{b}), we deduce that $u(t)/\Omega_{\omega}^0(t)\stackrel{t\to 0}{\longrightarrow }0$. Hence $u \in {\mathcal C}^{0}_{p, \omega}(0)$.

\smallskip
\noindent On the other hand, when $u\in {\mathcal C}^{0}_{p,\omega}(0)$,  then $u(t)\stackrel{t \to 0}{\longrightarrow}0$, because  $1/{\Omega}_{\omega}^0(t) \stackrel{t \to 0}{\longrightarrow} \infty$. Hence $u\in {\mathcal B}^{0}_{p,\omega}(0)={\mathcal A}^{0}_{p,\omega}(0)$.


\smallskip
\noindent {\bf iii):} \, Consider any sequence $t_n \searrow  0$.
 Then, for any fixed $t$, by using (\ref{est1}), we get the boundedness of $\{u(t_n)\}_n$. By Bolzano-Weierstrass Theorem, we can extract a converging subsequence,  which we will also denote by $\{u(t_n)\}_n$. Let $c$ be its limit.

\noindent By taking the limit as $n\to \infty$ in (\ref{est1}) we get
\begin{equation}\label{est2}
 \lvert u(t)-c  \rvert \leq \left(\int_0^t  {\lvert u^{'}(\tau) \rvert}^p\omega(\tau)d\tau\right)^{\frac{1}{p}}\cdot \Omega_{\omega}^0(t),
\end{equation}
which implies $u(t)\to c$ as $ t\to 0$.  Hence, any function $u \in D^{1,p}(\rp, \omega)$ has the limit as $t\to 0$, thus getting the  well-posedness of the trace operator $Tr^0(\cdot)$.

\noindent We have also proved that any function $u \in D^{1,p}(\rp, \omega)$ belongs to ${\mathcal B}^{0}_{p,\omega}(c)$ for some $c$.
This, together with (\ref{rdef}), gives the decomposition $ \displaystyle{D^{1,p}(\rp, \omega)=\bigcup_{c\in \mathbf{R}}{\mathcal R}^{0}_{p,\omega}(c)}$.

\smallskip
\noindent
{\bf iv):} \, Due to the existence of the limit of $u$ at zero, we can apply the estimate (\ref{czwartek}) with {$b:=a>0$ and $a:=0$,} thus getting
\begin{eqnarray}
 \lvert u(a)- Tr^0(u)  \rvert \le C_I \left( \int_0^\infty 
 { \lvert u^{'}(\tau)  \rvert}^p\omega d\tau\right)^{\frac{1}{p}}, {\rm where}\nonumber\\ C_I=\left( \int_{(0,a)} \omega^{-1/(p-1)}ds\right)^{1-\frac{1}{p}} .
\label{esttr}
\end{eqnarray}
Hence
\begin{eqnarray}
\|u\|_{D^{1,p}(\rp,\omega)}^{(a)}&=&\left( \int_0^\infty { \lvert u^{'}(\tau)  \rvert}^p\omega d\tau\right)^{\frac{1}{p}} +  \lvert u(a)  \rvert \label{equiv1}\\
&\le& \left( \int_0^\infty {  \lvert u^{'}(\tau)  \rvert}^p\omega d\tau\right)^{\frac{1}{p}} +  \lvert u(a)-Tr^0(u)  \rvert +  \lvert Tr^0(u)  \rvert \nonumber\\
&\stackrel{\eqref{esttr}}{\le}& \left( \int_0^\infty 
 {\lvert u^{'}(\tau)  \rvert}^p\omega d\tau\right)^{\frac{1}{p}} +C_I\left( \int_0^\infty {  \lvert u^{'}(\tau)  \rvert}^p\omega d\tau\right)^{\frac{1}{p}}+  \lvert Tr^0(u)  \rvert \nonumber\\
&=& (1+C_I) \left( \int_0^\infty  {  \lvert u^{'}(\tau)  \rvert}^p\omega d\tau\right)^{\frac{1}{p}}+  \lvert Tr^0(u)  \rvert \nonumber \\
&\le& (1+C_I) \, \|u\|_{D^{1,p}(\rp,\omega)}^{(0)} . \nonumber
\end{eqnarray}

\noindent On the other hand, by switching the rule of $a$ and $0$ {in \eqref{equiv1}, we obtain}
\begin{eqnarray*}
\|u \|_{D^{1,p}(\rp,\omega)}^{(0)}&=&\left( \int_0^\infty {  \lvert u^{'}(\tau)  \rvert}^p\omega d\tau\right)^{\frac{1}{p}} +  \lvert Tr^0(u)  \rvert \label{equiv2}\\
&\le& \left( \int_0^\infty   {\lvert u^{'}(\tau)  \rvert}^p\omega d\tau\right)^{\frac{1}{p}} +  \lvert Tr^0(u) -u(a)  \rvert +  \lvert u(a) \rvert \nonumber\\
&\le& \left( \int_0^\infty  {\lvert u^{'}(\tau)  \rvert}^p\omega d\tau\right)^{\frac{1}{p}} +C_I\left( \int_0^\infty  {\lvert u^{'}(\tau)
 \rvert}^p\omega d\tau\right)^{\frac{1}{p}}+  \lvert u(a)  \rvert\nonumber\\
&=& (1+C_I) \left( \int_0^\infty  {\lvert u^{'}(\tau)  \rvert}^p\omega d\tau\right)^{\frac{1}{p}}+ \lvert u(a)  \rvert \nonumber \\
&\le& (1+C_I) \, \|u\|_{D^{1,p}(\rp,\omega)}^{(a)}.\nonumber
\end{eqnarray*}

\smallskip
\noindent
This together with (\ref{equiv1}), yields the equivalence of all norms discussed, and completes the proof of the statement. \hfill$\Box$

\medskip
\noindent As a consequence of the above statement, we have is the following remarks.

\begin{rem}[representation of $\tilde{D}^{1,p}(\rp,\omega )$ for $\omega\in  B_p(0)$]\label{direct}\rm

\noindent Let $\omega: \rp \rightarrow \rp$, $\omega\in B_p(0)\cap C(\rp)$, $1<p<\infty$, and let {$Tr^{0}(\cdot)$} be as in (\ref{trzero}).
Then ${\mathcal R}^{0}_{p,\omega}(0)$ is a {Banach} subspace of $D^{1,p}(\rp,\omega )$ {(equipped with any of the norms $\|\cdot\|^{(a)}_{D^{1,p}(\rp,\omega )}$ where $a\in[0,\infty)$)}.
Moreover, the mapping
$$
D^{1,p}(\rp,\omega )\ni u \mapsto u-Tr^{0}(u)\in {\mathcal R}^{0}_{p,\omega}(0)
$$
is constant precisely on abstract classes in $\tilde{D}^{1,p}(\rp,\omega )$ (see Fact
\ref{norm-dirichlet}) and
 defines the isometric isomorphism between {($\tilde{D}^{1,p}(\rp,\omega ), \|\cdot\|^{*}_{\tilde{D}^{1,p}(\rp,\omega )})$} and \\
 $( {\mathcal R}^{0}_{p,\omega}(0), \|\cdot\|_{D^{1,p}(\rp,\omega)}^{(0)})$. In particular in  every  abstract class in $\tilde{D}^{1,p}(\rp,\omega )$,  there is the representative vanishing at zero and $\tilde{D}^{1,p}(\rp,\omega )$ {is} represented as
\begin{eqnarray*}
 \left \{ U= \left\{ \int_0^t v(\tau)d\tau + c\right\}_{  c\in\mathbf{R} }: v\in L^p(\rp,\omega)\, , \,
{ \| U\|_{\tilde{D}^{1,p}(\rp,\omega )}^{*} =\| v\|_{L^p(\rp,\omega)}} \right\} .
\end{eqnarray*}
 \end{rem}

\begin{rem}[asymptotic behaviour near zero]\label{assym}\rm
The statement $ii)$ in Theorem \ref{r-sets} and \eqref{est2} yield that if $\omega\in B_p(0)\cap C(\rp)$ {is  positive}, $1<p<\infty$,   then
for any  $u\in D^{1,p}(\rp, \omega)$
$$u(t)= Tr^0(u) + a(t)\cdot \Omega_\omega^0(t) \, ,$$

\noindent where $ a(t)$ is bounded, $a(t)\stackrel{t\to 0}{\to}0 $, and $\Omega_{\omega}^0 (\cdot )$ is as in (\ref{sladyzero}).
\end{rem}

\smallskip
\noindent
\subsection{Analysis in the case of $\omega \in B_p(\infty)$}

Let us assume that $ \omega \in B_p (\infty)\cap C(\rp)$ is  positive,  $1<p<\infty$. The aim of this section is to establish an analogous results to Theorem \ref{r-sets}, to representat the Dirichlet space through the trace operator, but in the case $ \omega \in B_p (\infty)$. The result stated below can be obatined by using very similar arguments to those used for the proof of Theorem \ref{r-sets}. Since we will deal with the $ B_p (\infty) $-condition,  we have to provide the analysis when $t$ is sufficiently large.
The proof is left to the reader with some general suggestions enclosed in order to treat this different setting:
\begin{itemize}
\item we first modify the appropriate definitions for the    sets from Theorem  \ref{r-sets};
\item in the proofs, we substitute the previously used limit conditions: $ t \searrow 0, t_n\searrow 0 $ by: $t \nearrow \infty$, $t_n\nearrow \infty$, respectively;
\item in place of \eqref{ajeden} we deal with the representation
\begin{equation*}\label{adwa}
u(t)=-\int_t^\infty u^{'}(\tau) d\tau ;
\end{equation*}
\item in place of  \eqref{est2} we deal with
\begin{equation}\label{est2infty}
 \lvert u(t)-c  \rvert \leq \left(\int_t^\infty  {\lvert u^{'}(\tau) \rvert}^p\omega(\tau)d\tau\right)^{\frac{1}{p}}\cdot \Omega_{\omega}^\infty (t),
\end{equation}
which forces $c=\displaystyle{\lim_{t\to\infty}u(t)}$.
\end{itemize}

\noindent
The following statement holds.

\smallskip
\begin{theo}[Asymptotic behaviour and trace  at infinity]\label{r-setsinfty}
Let $\omega :\mathbf{R}_+\rightarrow \mathbf{R}_+$,  $\omega\in C(\rp)\cap  B_p(\infty)$, $1<p<\infty$. For any $c\in \mathbf{R}$, let us consider the following subsets in
$D^{1,p}(\rp , \omega)$:
\begin{eqnarray}\label{sladyinfty}
{\mathcal R}^{\infty}_{p,\omega}(c)&:=&\left\{ u\in D^{1,p}(\rp , \omega),  u=  \int_t^{\infty} v(\tau)d\tau +c : v\in L^p(\rp,\omega)\right\}; \\
{\mathcal A}_{p,\omega }^{\infty}(c)&:=&\left\{ u\in D^{1,p}(\rp , \omega) :  \exists \, {t_n\nearrow \infty}: \lim_{n\to \infty} u(t_n)=c \right\};\nonumber\\
{\mathcal B}_{p,\omega }^{\infty}(c)&:=&  \left\{ u\in D^{1,p}(\rp , \omega ) :  \lim_{t\to \infty} u(t) =c\right\} ;\nonumber\\
{\mathcal C}^{\infty}_{p, \omega}(c)  &:=& \left\{ u\in D^{1,p}(\rp , \omega) :
 \lim_{t\to \infty} \frac{u(t)-c}{\Omega_{\omega}^{\infty} (t)} =0,\ \   {\rm sup}_{t>0} \frac{u(t)-c}{\Omega_{\omega}^{\infty} (t)}<\infty ,\right\}, \
  \nonumber
\end{eqnarray}

\noindent where $\displaystyle{{\Omega}_{\omega}^{\infty} (t)}=\left( \int_{t}^{\infty} \omega (\tau )^{\frac{-1}{p-1}  }d\tau\right)^{1-\frac{1}{p}}.$

\medskip
\noindent The following statements hold.
\begin{description}
\item[i)]  For any $c\in \mathbf{R}$, the set ${\mathcal R}^{\infty}_{p,\omega}(c)$ is a  closed subset in $D^{1,p}(\rp,\omega)$, equipped with any norm $\|\cdot\|_{D^{1,p}(\rp,\omega)}^{(a)}$ as in
\eqref{normaua}, where $a\in\rp$.
\item[ii)] For any $c\in \mathbf{R}$
\begin{eqnarray*}\label{rdefinfty}
{\mathcal A}_{p,\omega }^{\infty}(c)={\mathcal B}_{p,\omega }^{\infty}(c)={\mathcal C}^{}_{p, \omega}(c)= {\mathcal R}^{\infty}_{p,\omega}(c).
\end{eqnarray*}
\item[iii)]
 For every $u\in D^{1,p}(\rp,\omega)$ there is $c\in \mathbf{R}$ such that $u\in {\mathcal R}^{\infty}_{p,\omega}(c)$. In particular, the operator
\begin{equation}\label{trinfty}
Tr^{\infty}(u):= \lim_{t\to \infty} u(t)=: u(\infty),
\end{equation}
is well defined for every $u\in D^{1,p}(\rp, \omega)$ and $D^{1,p}(\rp ,\omega) = \cup_{c\in\mathbf{R}}  {\mathcal R}^{\infty}_{p,\omega}(c)$.

\noindent {Moreover,}
\begin{equation*}\label{dirrep}
D^{1,p}(\rp,\omega )  = \{ u= \int_t^{\infty} v(\tau)d\tau + c  : c\in\mathbf{R},  v\in L^p(\rp,\omega)\}.
\end{equation*}
\item[iv)] The quantity
$$ \| u\|_{D^{1,p}(\rp,\omega)}^{(\infty)} =
 \| u^{'}\|_{L^p(\rp,\omega)}+  \lvert Tr^{\infty}(u)  \rvert $$
is the norm on $D^{1,p}(\rp,\omega )$, which is equivalent to any norm $\| u\|^{(a)}$, where
$a\in \rp$.
\end{description}
\end{theo}

\begin{rem}[representation of $\tilde{D}^{1,p}(\rp,\omega )$ for $\omega\in  B_p(\infty)$]\label{directinfty}\rm
Let $\omega : \rp \rightarrow \rp$, $\omega\in B_p(\infty)\cap C(\rp)$ and  $1<p<\infty$,  and let $Tr^{\infty}(\cdot)$ be as in (\ref{trinfty}).
Then ${\mathcal R}^{\infty}_{p,\omega}(0)$ is a {Banach} subspace of $D^{1,p}(\rp,\omega )$   {(equipped with any of the norms $\|\cdot\|^{(a)}_{D^{1,p}(\rp,\omega )}$ where $a\in (0,\infty]$)}.
Moreover, the mapping
$$
D^{1,p}(\rp,\omega )\ni u \mapsto u-Tr^{\infty}(u)\in {\mathcal R}^{\infty}_{p,\omega}(0)
$$
is constant precisely on abstract classes in $\tilde{D}^{1,p}(\rp,\omega )$ (see Fact
\ref{norm-dirichlet}) and
 defines the isometric isomorphism between $\tilde{D}^{1,p}(\rp,\omega )$ and $( {\mathcal R}^{\infty}_{p,\omega}(0), \|\cdot\|_{D^{1,p}(\rp,\omega)}^{(\infty)})$. In particular in  every  abstract class in $\tilde{D}^{1,p}(\rp,\omega )$ there is the representative vanishing at infinity and $\tilde{D}^{1,p}(\rp,\omega )$ {represents as}
\begin{eqnarray*}
 \left \{ U= \left\{ \int_t^\infty v(\tau)d\tau + c\right\}_{  c\in\mathbf{R} }: v\in L^p(\rp,\omega)\, , \,
 \| U\|^{{*}}_{\tilde{D}^{1,p}(\rp,\omega )} {=}\| v\|_{L^p(\rp,\omega)} \right\} .
\end{eqnarray*}
 \end{rem}

 \begin{rem}[asymptotic behaviour near infinity]\label{assyminfty}\rm
The statement ii) in Theorem \ref{r-setsinfty} and \eqref{est2infty} yield that if $\omega : \rp \rightarrow \rp$, $\omega\in B_p(\infty)\cap C(\rp)$  and  $1<p<\infty$, then
for {every}  $u\in D^{1,p}(\rp, \omega)$
$$u(t)= Tr^\infty(u) + a(t)\cdot \Omega_\omega^\infty(t) \, ,$$
where $ a(t)$ is bounded, $a(t)\stackrel{t\to \infty}{\to}0 $, and $\Omega_{\omega}^{\infty} (\cdot )$ is as in (\ref{sladyinfty}).
\end{rem}

\section{Characterization  of ${D}^{1,p}_0(\rp ,\omega)$ and density results}\label{densityres}

\subsection{The space  ${D}^{1,p}_0(\rp ,\omega)$ and first density results}

We start with the following definition.

\begin{defi}[the space ${D}^{1,p}_0(\rp ,\omega)$, the case of $\omega\in B_p$]\label{dzero-general}\rm~\\When $\omega\in B_p$,  $1<p<\infty$,
by ${D}^{1,p}_0(\rp ,\omega)$ we will denote
the subset of all functions $u \in D^{1,p}(\rp ,\omega)$, for which there exists the sequence
$\{ \phi_n\}\subseteq C_0^\infty (\rp)$, which satisfies:
 $$\phi_n^{'}\rightarrow u^{'}\ \ {\rm  in}\ \
L^p(\rp,\omega) \ {\rm and } \ \phi_n\to u \ \hbox{\rm uniformly on compact sets in}\ \rp .$$
\end{defi}
\noindent By { Fact \ref{bp-dirichlet}},
 ${D}^{1,p}_0(\rp ,\omega)$  is the same as the completion of $C_0^\infty(\rp)$ in the space $D^{1,p}(\rp ,\omega)$ equipped with any of the norms $\|\cdot\|_{D^{1,p}(\rp,\omega)}^{(a)}$, where $a\in\rp$ can be taken arbitrary.
In particular, it is the Banach subspace of  \\$(D^{1,p}(\rp,\omega), \|\cdot\|_{D^{1,p}(\rp,\omega)}^{(a)})$, with an arbitrary $a\in\rp$.

\bigskip
\noindent
The following fact is rather obvious to the specialists, but for reader's convenience we submit its proof.

\begin{lem} \label{gestosci}
Let $\omega: \rp \rightarrow \rp$, $\omega\in C(\rp)$ and $1< p<\infty$.
Then for any $a\in \rp$
\begin{eqnarray*}\label{jed-a}
D_0
^{1,p}(\rp , \omega)=\overline{Lip_c}^{\|\cdot\|^{(a)}_{D^{1,p}(\rp , \omega)}}= \overline{\left( D^{1,p}(\rp , \omega)\right)_{c}}^{\|\cdot\|^{(a)}_{D^{1,p}(\rp , \omega)}}.
\end{eqnarray*}
\end{lem}

\noindent
{\bf Proof.}
Clearly,
$$
C_0^\infty (\rp) \subseteq Lip_c (\rp ) \subseteq \left(D^{1,p}(\rp ,\omega)\right)_{c}.
$$
Hence,  it suffices to show that $\overline{\left( D^{1,p}(\rp , \omega)\right)_{c}}^{\|\cdot\|^{(a)}_{D^{1,p}(\rp , \omega)}} \subseteq D^{1,p}_{0}(\rp ,\omega)$. For that, take $u\in (D^{1,p}(\rp ,\omega))_{c}$ with the support $[a,b]\subseteq \rp$.
As on compactly supported sets $\omega \sim 1$,
therefore $u\in D^{1,p}(\rp)$ and $u$ is compactly supported. By standard convolution arguments, the convolutions
$u_\epsilon (x):= \phi_\epsilon *u$, with the classical mollifier functions $\phi_\epsilon (x)=\epsilon^{-1}\phi (x/\epsilon)$,  where
$\phi\in C_0^\infty (\mathbf{R})$, $0\le \phi\le 1$, ${\rm supp}u\subseteq [-1,1]$ and $\int_{\mathbf{R}} \phi dx =1$, converge to $u$ in the topology of $D^{1,p}(\rp ,\omega)$.
Moreover,
their supports are subsets of $J:=[a/2, 3/2b]$ for the sufficiently small $\epsilon$'s. Again, as $\omega\sim 1$  on $J$, therefore
$u_\epsilon$'s converge to $u$ also in $D^{1,p}(\rp ,\omega)$.
This shows that $u\in D^{1,p}_{0}(\rp ,\omega)$.
\hfill$\Box$

\bigskip
\noindent
In the preceding sections we will analyze independently the  cases:
$\omega\in B_p(0)$ and $\omega\in B_p(\infty)$.

\subsection{The case of $\omega \in B_p(0)$} \label{sec-density-0}

Let $\omega : \rp \rightarrow \rp$, $\omega\in B_p(0)\cap C(\rp)$ and  $1<p<\infty$, and let us consider
the space $D^{1,p}_0(\rp, \omega)$ as in Definition \ref{dzero-general}. According to Theorem \ref{direct}, part iv), it is the completion of
$C_0^\infty (\rp)$ in $D^{1,p}(\rp, \omega)$ equipped with the norm
$\|\cdot\|^{(0)}$. As $C_0^\infty (\rp)\subseteq \mathcal{R}^{0}_{p,\omega}(0)$ and by Theorem \ref{r-sets} $\mathcal{R}^{0}_{p,\omega}(0)$ is a closed subspace in $D^{1,p}(\rp, \omega)$, we deduce that
 \begin{equation}\label{a2}
 D_0^{1,p}(\rp,\omega)\subseteq \mathcal{R}^0_{p,\omega}(0),\ {\rm when}\
 \omega\in B_p(0)\cap C(\rp)\ \hbox{\rm is  positive}.
  \end{equation}

\medskip
\noindent
We address the question about density:
$$  {\rm When}\ \ \ D_0^{1,p}(\rp, \omega) = {\mathcal{R}}^{0}_{p,\omega}(0) \, \, ?$$
 The statement given below answers on this question.

\bigskip

\begin{theo}[characterization of weights for $D^{1,p}_0(\rp,\omega)={\mathcal R}_{p,\omega}^{0}(0)$]\label{characterization}~\\
Let $\omega:\rp\rightarrow \rp, \omega  \in B_p(0)\cap C(\rp),  1<p<\infty$.
Then
\begin{equation*}
D_0^{1,p}(\rp,\omega )= {\mathcal R}_{p,\omega}^{0}(0) \Longleftrightarrow \omega\not\in B_p(\infty).
\end{equation*}
\end{theo}

\smallskip
\noindent
The proof will be based on the following lemma, whose proof is submitted in the Appendix for reader's convenience.

\begin{lem}[Energy minimizer of nontrivial constraint at left end]\label{lemat1}~\\ Let $\omega : \rp \rightarrow \rp$, $\omega \in  C(\rp)$, $0<k<K<\infty, 0\neq a\in \mathbf{R}$, $1<p<\infty$, and consider energy functional
\begin{eqnarray}\label{energy}
E_\omega (\phi) :=  \int_k^K  {  \lvert \phi^{'}(t)  \rvert}^p\omega(t)dt,\ \ \phi \in W^{1,p}((k,K)), \ \phi (k)=a, \phi (K)=0.
\end{eqnarray}
Then the minimum of $E_\omega (\cdot)$ is achieved at
\begin{equation*}
\phi_{(k,K,a)}(t):=a\left(\int_k^K \omega^{\frac{-1}{p-1}}(\tau)d\tau\right)^{-1} \int_t^K \omega^{\frac{-1}{p-1}}(\tau)d\tau .
\end{equation*}
\end{lem}

\noindent
We are now to prove Theorem \ref{characterization}.

\smallskip
\noindent
{\bf Proof of Theorem \ref{characterization}:}
''$\Longleftarrow$''  ($E_\omega$- caloric approximation):\\
We have called this part of the proof ``$\omega$- caloric approximation'', because  the construction of the approximation sequence involves the energy minimizers of (\ref{energy}).\\
Assume that $\omega\in C(\rp)\cap B_p(\infty)\setminus B_p(0)$. Thanks to \eqref{a2} we only have to prove that $\mathcal{R}^0_{p,\omega}(0)\subseteq D^{1,p}_0(\rp,\omega)$.

\noindent The proof follows by two steps.

\smallskip
\noindent
{\sc Step 1. Reduction argument.}
We show that it suffices to prove that any $u\in D^{1,p}(\rp,\omega)$, such that
\begin{equation}\label{spacial}
u\equiv 0\ \ \hbox{\rm  near zero and}\ \ u\equiv 1\ \
{\rm on}\  (k,\infty),\ \hbox{\rm for some}\  k>0,
\end{equation}
belongs to $D_0^{1,p}(\rp,\omega)$.

Indeed, let us take
$u\in \mathcal{R}^{0}_{p,\omega}(0)$, we notice that functions in the form
\begin{equation}\label{zbieznoscbpzero}
\tilde{u}_n (t):= \int_0^t u^{'}(\tau)\chi_{(\frac{1}{n},n)} (\tau)d\tau
\end{equation}
converge to $u$ in $(D^{1,p}(\rp,\omega), \|\cdot\|^{(0)})$ (see
 \eqref{normazero}),  they are zero near zero and constant near infinity.
 Clearly, if that constant equals zero,  according to Lemma \ref{gestosci}, we have $\tilde{u}_n\in D_0^{1,p}(\rp,\omega)$. In the other case, we are left with the proof that $\tilde{u}_n\in D_0^{1,p}(\rp,\omega)$. Obviously, it suffices to consider  $C\tilde{u}_n$ instead of $\tilde{u}_n$, with  constant $C$ such that $C\tilde{u}_n\equiv 1$ near infinity.

\noindent
{\sc Step 2. Proof in the special case.}
We prove that any $u\in D^{1,p}(\rp,\omega)$ as in (\ref{spacial})
belongs to $D_0^{1,p}(\rp,\omega)$.

\smallskip
\noindent
Let  $u\in D^{1,p}(\rp,\omega)$ satisfy (\ref{spacial}). For any $n\in\mathbf{N}$ and $k<t_n$, let
\begin{eqnarray*}
u_n(t):= \left\{
\begin{array}{ccc}
u(t) &{\rm if}& t<k\\
\phi_{(k,t_n,1)}(t) & {\rm if}& t\in [k,t_n]\\
0 & {\rm if} & t>t_n
\end{array}
\right.
\end{eqnarray*}
where $\phi_{(k,t_n,1)}$ is as in Lemma \ref{lemat1} and $t_n\nearrow \infty$. Clearly, the $u_n$'s are compactly supported.
We will show that
\begin{equation}\label{zbieznosc-zwarte-nosniki}
u_n\stackrel{n\to\infty}{\rightarrow} u \ {\rm in}\  \ (D^{1,p}(\rp,\omega), \|\cdot\|_{D^{1,p}(\rp,\omega)}^{(0)}),
\end{equation}
 which, together with Lemma \ref{gestosci},  will close the assertion for this part of the statement.
\noindent
 We  have:
\begin{eqnarray}\label{wybortn}
\int_0^\infty  {\lvert (u_n-u)^{'}(t)  \rvert}^p \omega (t)dt =
\int_k^{t_n}{  \lvert (1-\phi_{(k,t_n,1)}(t))^{'}  \rvert}^p \omega (t) dt =
\nonumber\\
=\int_k^{t_n}{  \lvert \phi_{(k,t_n,1)}(t)^{'}  \rvert}^p \omega (t) dt =
\int_k^{t_n} \frac{\omega^{\frac{-p}{p-1}}(t)}{\left(
\int_k^{t_n}\omega^{\frac{-1}{p-1}}(\tau ) d\tau
\right)^p}
\cdot \omega(t)dt\\
=
\left( \int_k^{t_n} {  \lvert \omega^{\frac{-1}{p-1}}(\tau)  \rvert}^p d\tau\right)^{-p}
\int_k^{t_n} \omega^{\frac{-1}{p-1}}(\tau)d\tau =
\left( \int_k^{t_n}\omega^{\frac{-1}{p-1}}(\tau) d\tau\right)^{1-p}.\nonumber
\end{eqnarray}
As $\omega\not\in B_p(\infty)$ and $p>1$,
$$
\int_k^{t_n}\omega^{\frac{-1}{p-1}}(\tau) d\tau \stackrel{n\to\infty}{\rightarrow }\infty ,\ \ \hbox{\rm consequently}\ \ \left( \int_k^{t_n}\omega^{\frac{-1}{p-1}}(\tau) d\tau\right)^{1-p}\stackrel{n\to\infty}{\rightarrow } 0.
$$
{This implies \eqref{zbieznosc-zwarte-nosniki}.}

\bigskip
\noindent
``$\Longrightarrow$'':\\
Suppose that $D_0^{1,p}(\rp,\omega )= {\mathcal R}^{1,p}(0)$. We will show that $\omega\not\in B_p(\infty)$.

\noindent
Clearly, the function
\begin{eqnarray*}
u(t):=\left\{
\begin{array}{ccc}
0 &{\rm if}& t\le 1\\
t-1 &{\rm if} & t\in [1,2]\\
1 &{\rm if}& t>2
\end{array}
\right.
\end{eqnarray*}
belongs to $\mathcal{R}^{0}_{p,\omega}(0) (= D_0^{1,p}(\rp, \omega))$, and so there is the sequence
$\{ u_n\}_{n\in\mathbf{N}}\subseteq C_0^\infty (\rp)$ such that
$u_n\rightarrow u$ in $(D^{1,p}(\rp,\omega), \| \cdot\|_{D^{1,p}(\rp,\omega)}^{(2)})$, and in all the equivalent norms.
In particular
$$
\xi_n:= u_n(2)\stackrel{n\to\infty}{\rightarrow} 1 \ {\rm and}\
E_\omega (u_n):= \int_{(2,t_n)} {\lvert u_n^{'}(\tau) \rvert}^p \omega(\tau) d\tau \stackrel{n\to\infty}{\rightarrow} 0 ,
$$
for any  $t_n>2$.
Let  $\{ t_n\}_{n\in\mathbf{N}}$ be any sequence such that
${\rm supp}\, u_n\subseteq (0,t_n)$ and $t_n\nearrow\infty$ as $n\to\infty$.
According to Lemma \ref{lemat1},  the energy $E_\omega (\phi_{(2,t_n,\xi_n)})$ cannot be larger than $E_\omega (u_n)$. Therefore we also have
$$
I_n:=E_\omega (\phi_{(2,t_n,\xi_n)})= \int_{(2, t_n)} { \lvert (\phi_{(2,t_n,\xi_n)})^{'}(\tau) \rvert }^p \omega(\tau) d\tau \stackrel{n\to\infty}{\rightarrow} 0 ,
$$
Since
\begin{eqnarray*}
I_n = {\lvert  \xi_n \rvert }^p \int_{(2,t_n)}  \left(\frac{ \omega^{-\frac{1}{p-1}}}{\int_{(2,t_n)} \omega^{-\frac{1}{p-1}}(\tau)d\tau}\right)^p\omega(t)dt =  {\lvert \xi_n \rvert }^p
\left( \int_{(2,t_n)}  \omega^{-\frac{1}{p-1}}(\tau)d\tau\right)^{1-p},
\end{eqnarray*}
it follows that the above converges to zero if and only if
$$
\lim_{t\to\infty} \int_{(2,t)} \omega^{-\frac{1}{p-1}}(\tau)d\tau =\infty,
$$
equivalently $\omega\not\in B_p(\infty)$.
This completes the proof of the statement.
\hfill$\Box$.

\subsection{Analysis in the case of $\omega \in B_p(\infty)$} \label{sec-density-infty}


Let us assume that $\omega : \rp \rightarrow \rp$, $\omega\in B_p(\infty)\cap C(\rp)$ and $1<p<\infty$. Our aim is to analyze the properties of the space
$D^{1,p}_0(\rp, \omega)$, the completion of
$C_0^\infty (\rp)$ in $D^{1,p}(\rp, \omega)$ equipped with the norm
$\|\cdot\|_{D^{1,p}(\rp,\omega)}^{(\infty)}$, defined  in Theorem \ref{r-setsinfty}. As $C_0^\infty (\rp)\subseteq \mathcal{R}^{\infty}_{p,\omega}(0)$, where the latter space, due to Theorem \ref{r-setsinfty},  is a closed subspace in $D^{1,p}(\rp, \omega)$, we deduce that
 \begin{equation}\label{a2infty}
 D_0^{1,p}(\rp,\omega)\subseteq \mathcal{R}^\infty_{p,\omega}(0).
  \end{equation}

\smallskip
\noindent
The goal of this section is the following characterization theorem.

\begin{theo}[characterization of weights for $D^{1,p}_0(\rp,\omega)={\mathcal R}_{p,\omega}^{\infty}(0)$]\label{characterizationinfty}~\\
Let $\omega:{\rp \rightarrow \rp,} \, \omega \in B_p(\infty)\cap C(\rp)$ and $1<p<\infty$.
Then
\begin{equation*}
D_0^{1,p}(\rp,\omega )= {\mathcal R}_{p,\omega}^{\infty}(0) \Longleftrightarrow \omega\not\in B_p(0).
\end{equation*}
\end{theo}

\noindent {The proof is an easy modification of the proof of Theorem \ref{characterization}, so we only sketch it. We}  start by stating the following  result similar to that of Lemma \ref{lemat1}. Its proof can be easily obtained either by the internal symmetry argument, or by suitable modification of the proof of Lemma \ref{lemat1}. {We leave it to the reader.}

\begin{lem}[Energy minimizer of nontrivial constraint at right end]\label{lemat2}~\\ Let $\omega:{\rp \rightarrow \rp ,}  \, \omega  \in B_p\cap C(\rp)$ and $1<p<\infty$, $0<k<K<\infty, 0\neq a\in \mathbf{R}$, and consider energy functional
\begin{eqnarray*}\label{energy1}
\tilde{E}_\omega (\phi) :=  \int_k^K { \lvert \phi^{'}(t) \rvert }^p\omega(t)dt,\ \ \phi \in W^{1,p}((k,K)), \ \phi (k)=0, \phi (K)=a.
\end{eqnarray*}
Then the minimum of $\tilde{E}_\omega (\cdot)$ is achieved at
\begin{equation*}
\tilde{\phi}_{(k,K,a)}(t):=a\left(\int_k^K \omega^{\frac{-1}{p-1}}(\tau)d\tau\right)^{-1} \int_k^t \omega^{\frac{-1}{p-1}}(\tau)d\tau .
\end{equation*}
\end{lem}

\bigskip
\noindent
We are in position to sketch the proof of Theorem \ref{characterizationinfty}.

\bigskip
\noindent
{\bf Proof of Theorem \ref{characterizationinfty}.} ``$\Longleftarrow$'' ($\tilde{E}_\omega$ - caloric approximation):

\noindent
Assume that  $\omega\in B_p(\infty)\setminus B_p(0)$. As \eqref{a2infty} holds, we only have to prove that $\mathcal{R}^\infty_{p,\omega}(0)\subseteq
D^{1,p}_0(\rp, \omega)$.

Let $u\in {\mathcal R}_{p,\omega}^{\infty}(0)$. We will show that
$u\in D_0^{1,p}(\rp,\omega)$.
The proof consists of two steps.

\smallskip
\noindent
{\sc Step 1. Reduction argument.}
We show that it suffices prove that any $u\in D^{1,p}(\rp,\omega)$, such that
\begin{equation}\label{spacialinfty}
u\equiv 1\ \ \hbox{\rm  on some} \ (0,c)\ {\rm where}\ c>0,\ \ \ u\equiv 0\ \
\hbox{\rm  near}\ \infty ,
\end{equation}
belongs to $D_0^{1,p}(\rp,\omega)$.

\noindent To this aim, we note that functions
\begin{equation*}\label{zbieznoscbpzeroinfty}
\tilde{u}_n (t):= - \int_t^\infty u^{'}(\tau)\chi_{(\frac{1}{n},n)} (\tau)d\tau
\end{equation*}
are proportional to functions as in \eqref{spacialinfty}, and they
converge to $u$ in $(D^{1,p}(\rp,\omega), \|\cdot\|^{(\infty)})$.

\smallskip
\noindent
{\sc Step 2. Proof in the special case.}\\
{We} prove that any $u\in D^{1,p}(\rp,\omega)$ as in (\ref{spacialinfty})
belongs to $D_0^{1,p}(\rp,\omega)$. To this purpose, let us consider the following sequence dealing with $0<s_n<c$, $s_n\searrow 0$:

\begin{eqnarray*}
u_n(t):= \left\{
\begin{array}{ccc}
u(t) &{\rm if}& t>c\\
\tilde{\phi}_{(s_n,c,1)}(t) & {\rm if}& t\in [s_n,c]\\
0 & {\rm if} & t<s_n,
\end{array}
\right.
\end{eqnarray*}
where  $\tilde{\phi}_{(s_n,c,1)}$ is as in Lemma \ref{lemat2}. Obviously the $u_n$'s are compactly supported and so they belong to $D_0^{1,p}(\rp, \omega)$, by Lemma \ref{gestosci}.
By similar computations as in (\ref{wybortn}), we get
$$
\int_0^\infty { \lvert (u_n-u)^{'}(t) \rvert}^p \omega (t)dt =
\left( \int_{s_n}^c \omega^{\frac{-1}{p-1}}(\tau) d\tau\right)^{1-p}\stackrel{n\to\infty}{\rightarrow}0,
$$
because $\omega\not\in B_p(0)$. Therefore  $u \in D_0^{1,p}(\rp,\omega)$.

\smallskip
\noindent
``$\Longrightarrow $'':\\
We prove that the condition $D_0^{1,p}(\rp,\omega )= {\mathcal R}^{1,p}(\infty)$ forces the condition $\omega\not\in B_p(0)$.
To this aim, assume that $D_0^{1,p}(\rp,\omega )= {\mathcal R}^{1,p}(\infty)$ and let
\begin{eqnarray*}
u(t):=\left\{
\begin{array}{ccc}
1 &{\rm if}& t\le 1\\
2-t &{\rm if} & t\in [1,2]\\
0 &{\rm if}& t>2
\end{array}
\right.
\end{eqnarray*}
Then $u \in {\mathcal{R}}_{p,\omega}^{\infty}(0)\subseteq{D_0^{1,p}(\rp,\omega)} $. Hence, there exists a sequence
$\{ u_n\}_{n\in\mathbf{N}}\subseteq C_0^\infty (\rp)$ such that
$u_n\rightarrow u$ in $(D^{1,p}(\rp,\omega), \| \cdot\|_{D^{1,p}(\rp,\omega)}^{(1)})$.
In particular
$$
\xi_n:= u_n(1)\stackrel{n\to\infty}{\rightarrow} 1 \ {\rm and}\
\tilde{E}_\omega (u_n):= \int_{s_n}^1 {\lvert  u_n^{'}(\tau) \rvert }^p \omega(\tau) d\tau \stackrel{n\to\infty}{\rightarrow} 0 ,
$$
where ${\rm supp}\, u_n\subseteq (s_n,\infty)$, for some  $s_n\searrow 0$.
By Lemma \ref{lemat2}, the energy $I_n:=\tilde{E}_\omega (\tilde{\phi}_{(s_n,1,\xi_n)})$ cannot be larger than $\tilde{E}_\omega (u_n)$, therefore
$$
I_n :=  {\lvert  \xi_n \rvert }^p
\left( \int_{s_n}^1  \omega^{\frac{-1}{p-1}}(\tau)d\tau\right)^{1-p} \stackrel{n\to\infty}{\rightarrow} 0.
$$
{Consequently} $\omega\not\in B_p(0)$, which completes the proof.
\hfill$\Box$

\subsection{Analytic description of $D^{1,p}_0(\rp,\omega)$ in general case}

\noindent
Our main statement in this section reads as follows.

\begin{theo}[description  of $D^{1,p}_0(\rp,\omega)$ for all admitted weights]
\label{glowne}
Let $\omega : \rp \rightarrow \rp$, $\omega\in  C(\rp)$ and $1<p<\infty$.
Then we have.
\begin{description}
\item[i)]
If $\omega\in B_p(0)\setminus B_p(\infty)$, then
$$ D^{1,p}_0(\rp,\omega)= \mathcal{R}^0_{p,\omega}(0).$$
\item[ii)]
If $\omega\in B_p(\infty)\setminus B_p(0)$,  then
$$ D^{1,p}_0(\rp,\omega)= \mathcal{R}^{\infty}_{p,\omega}(0).$$
\item[iii)]
If $\omega\not\in B_p(0)\cup B_p(\infty)$,  then
$$ D^{1,p}_0(\rp,\omega)= D^{1,p}(\rp,\omega).$$
\item[iv)]
If $\omega \in B_p(0)\cap B_p(\infty)$,  then
$$ D^{1,p}_0(\rp,\omega)= \mathcal{R}^0_{p,\omega}(0)\cap \mathcal{R}^\infty_{p,\omega}(0).$$
\end{description}
\end{theo}
\noindent
{\bf Proof.}
{{\bf i) and ii):}}
Statements i) and ii) have been already obtained in Theorems \ref{characterization} and \ref{characterizationinfty}. We are left with the proofs of parts iii) and iv).

\smallskip
\noindent
{\bf iii):} Assume that $\omega\not\in B_p(0)\cup B_p(\infty)$.
 Let $u\in D^{1,p}(\rp,\omega)$ and
let us consider the Lipschitz resolution of the unity on $\rp$:  $\phi_0,\phi_1$, defined by
$$\phi_1 (t) =1\chi_{(0,1)}+ (-t+2)\chi_{[1,2]},\ \  \phi_0(t):=1-\phi_1.$$ We have to prove that $u\in  D^{1,p}_0(\rp,\omega)$. As  $u=\phi_0u+\phi_1u$, it suffices to consider the following cases:
a)  $u\equiv 0$ near $0$  and b)
$u\equiv 0$ near $\infty$. \\
In case a), suppose that $u\equiv 0$ on $(0,a]$ for some $a>0$. Then functions as in (\ref{zbieznoscbpzero}) converge to $u$ in $(D^{1,p}(\rp,\omega),\| \cdot\|_{D^{1,p}(\rp,\omega)}^{(a)})$,
they { are zero} on $(0,a]$ and constant when $t>n$. Therefore the proof reduces to the case of $u\equiv 0$ near zero and $u\equiv Const$ near $\infty$. Then we repeat all the arguments
from the proof of Theorem \ref{characterization}, Step 2, in the proof of the implication ``$\Longleftarrow$''.
\\
In case b), the argument at the beginning of the proof of Theorem
\ref{characterizationinfty} reduces that case to the situation when
$u\equiv Const$ near $0$ and $u\equiv 0$ near to infinity. In that case we use precisely the same arguments as in the proof of Theorem \ref{characterizationinfty}, Step 2 in  part ``$\Longleftarrow$''.

\smallskip
\noindent
{\bf iv): } Let $\omega\in B_p(0)\cap B_p(\infty)$. 

\noindent
Then we have
$D_0^{1,p}(\rp,\omega)\subseteq \mathcal{R}^{0}_{p,\omega}(0)\cap
\mathcal{R}^{\infty}_{p,\omega}(0)$, by \eqref{a2} and \eqref{a2infty}.
Therefore it suffices to show that the converse inclusion holds. For that, assume that $u\in \mathcal{R}^{0}_{p,\omega}(0)\cap \mathcal{R}^{\infty}_{p,\omega}(0)$.
Note that in particular $u$ is bounded and, by Fact \ref{zanurzenia},
$u^{'}\in L^1_{loc}([0,\infty))\cap L^1_{loc}((0,\infty])=L^1((0,\infty))$. Moreover,
$$
u(t)=\int_0^t u^{'}(\tau )d\tau = - \int_t^\infty u^{'}(\tau )d\tau .
$$
Thus we have
$$
\int_{(0,\infty )} u^{'}(\tau )d\tau =0 \ {\rm and}\ u^{'}\in L^1 (\rp) .
$$
For each $n\in\mathbf{N}$, consider
$$
u_n (t):= \int_0^t \chi_{(\frac{1}{n},n)} \left( u^{'}(\tau ) - c_n\right) d\tau, \ \ {\rm where}\
c_n =\frac{1}{n-\frac{1}{n}} \int_{\frac{1}{n}}^n u^{'}(\tau )d\tau .
$$
Then $$u_n^{'}(\tau) =\chi_{(\frac{1}{n},n)} \left( u^{'}(\tau ) - c_n\right),\ \
\int_0^{\infty} u_n^{'}(\tau )d\tau =0,\  \ u_n^{'}\in L^1 (\rp)\cap L^p(\rp,\omega).$$  {Consequently}
$$
u_n(t)= \int_0^t u_n^{'}(\tau )d\tau = - \int_t^\infty u_n^{'}(\tau )d\tau ,\ \ u_n^{'}\in L^p(\rp,\omega).
$$
Moreover, we have $u_n\stackrel{n\to\infty}{\rightarrow} u$ in $D^{1,p}(\rp,\omega)$ and $u_n$ is supported in $\displaystyle{[\frac{1}{n},n]}$.
This together with Lemma \ref{characterization}
 implies that $u\in D_0^{1,p}(\rp,\omega)$ and ends the proof of the statement.
\hfill$\Box$

\subsection{Sharpness in Theorems:  \ref{r-sets},  \ref{r-setsinfty} and \ref{glowne}}

This section is devoted to prove sharpness  in statements $i)$ in Theorems \ref{r-sets} and \ref{r-setsinfty},
and the converse implications in Theorem \ref{glowne}.

\smallskip
\noindent We have the following result.

\noindent
\begin{theo}[sharpness in Theorems \ref{r-sets},  \ref{r-setsinfty}, parts i)] \label{domknietosczbiorowzero}~\\
{Let  $\omega : \rp \rightarrow \rp$, $\omega\in C(\rp)$ and $1<p<\infty$. Then the following statements hold.}
\item[\bf {\rm i})]
 $\mathcal{R}_{p,\omega}^{0}(0)$ is a closed subset in $D^{1,p}(\rp,\omega)$ $\Longleftrightarrow$
$\omega\in B_p(0)$.
\item[\bf {\rm ii})]  $\mathcal{R}_{p,\omega}^{\infty}(0)$ is a closed subset in $D^{1,p}(\rp,\omega)$ $\Longleftrightarrow$
$\omega\in B_p(\infty)$.
\end{theo}
\noindent
{\bf Proof.} The implications ``$\Longleftarrow$'' were already proved in Theorem \ref{r-sets},  \ref{r-setsinfty}, parts i). We are left with the proofs of the  converse implications.
\\
``$\Longrightarrow :$''

\smallskip
\noindent {\bf i){:}}
We argue by  contradiction. Assume that the implication does not hold, that is  $\mathcal{R}_{p,\omega}^{0}(0)$ is a closed subset in $D^{1,p}(\rp,\omega)$,  but $\omega\not\in B_p(0)$.
We have either  a) $\omega \in B_p(\infty)$ \, or\,  b) $\omega \not\in B_p(\infty)$.

\noindent If  $\rm a)$ holds {, then by part ii) of Theorem \ref{glowne},}
 $D^{1,p}_0(\rp, \omega) = \mathcal{R}^{\infty}_{p,\omega}(0)$.

\noindent The function
\begin{eqnarray}\label{exfunc1}
u(x):= \left\{
\begin{array}{ccc}
1 & {\rm if} & x\in (0,1)\\
-x +2  & {\rm if} & x\in [1,2)\\
0 & {\rm if} & x\ge 2
\end{array}
\right.
\end{eqnarray}
belongs to $\mathcal{R}^{\infty}_{p,\omega}(0)$,  and consequently to $D^{1,p}_0(\rp, \omega)$. Hence,  it can be approximated in $D^{1,p}(\rp, \omega)$ by functions which are zero near zero, that is belonging to $\mathcal{R}^{0}_{p,\omega}(0)$. At the same time  their limit $u(x)\not\in \mathcal{R}^{0}_{p,\omega}(0)$. Therefore $\mathcal{R}^{0}_{p,\omega}(0)$ cannot be closed. The  condition $\rm a)$ is not possible.

\noindent Let us suppose that  $\rm b)$ holds.
By Theorem \ref{glowne}, part iii), we have $D^{1,p}_0(\rp, \omega) = D^{1,p}(\rp, \omega)$. Thus, the function $u$ from (\ref{exfunc1}) belongs to $D_0^{1,p}(\rp, \omega)$. We can argue as in the previous case,  obtaining the contradiction that $\mathcal{R}^0_{p,\omega}$ is closed.
Therefore necessarily  $\omega \in B_p(0)$.


\medskip
\noindent {\bf ii):}
Let us assume that the implication ``$\Longrightarrow$''
 does not hold, that is
$\mathcal{R}_{p,\omega}^{\infty}(0)$ is a closed subset in $D^{1,p}(\rp,\omega)$ but $\omega\not\in B_p(\infty)$.

\noindent We have either a) $\omega \in B_p(0)$ \, or\,  b) $\omega \not\in B_p(0)$. {Both conditions,} by Theorem \ref{glowne},  imply that $D^{1,p}_0(\rp, \omega) = \mathcal{R}^{0}_{p,\omega}(0)$ or $D^{1,p}_0(\rp, \omega) = D^{1,p}(\rp, \omega)$, respectively.

\noindent Let us consider the function

{
\begin{eqnarray*}\label{exfunc2}
u(x):= \left\{
\begin{array}{ccc}
0 & {\rm if} & x\in (0,1)\\
x -1  & {\rm if} & x\in [1,2)\\
1 & {\rm  if} & x\ge 2
\end{array}
\right.
\end{eqnarray*}
}
It belongs to $\in D_0^{1,p}(\rp,\omega)$.
By arguments as in the proof of part  i), we get a contradiction in both cases: $\rm a)$ and $\rm b)$, which proves $\rm ii)$.

\noindent
The proof of the statement is {complete}.
\hfill$\Box$

\medskip
\noindent Let us proceed by proving the {converse} implications  in Theorem \ref{glowne}.

\noindent We state the following

\begin{theo}[sharpness of conditions on weights in Theorem \ref{glowne}]\label{zamkniecieimplikacji}~\\
 Let $\omega : \rp \rightarrow \rp$, $\omega\in C(\rp)$ and  $1<p<\infty$. Then we have.
\noindent \item[\bf i)]
 $D_0^{1,p}(\rp,\omega) = \mathcal{R}_{p,\omega}^{0}(0) \, \Longleftrightarrow \,
\omega\in B_p(0)\setminus B_p(\infty)$.
\item[\bf ii)]
$D^{1,p}_0(\rp,\omega)= \mathcal{R}^{\infty}_{p,\omega}(0)\, \Longleftrightarrow \, \omega\in B_p(\infty)\setminus B_p(0)$.
\item[\bf iii)]
 $ D^{1,p}_0(\rp,\omega)= D^{1,p}(\rp,\omega) \, \Longleftrightarrow \, \omega\not\in B_p(0)\cup B_p(\infty)$.
\item[\bf iv)]
 $ D^{1,p}_0(\rp,\omega)= \mathcal{R}^0_{p,\omega}(0)\cap \mathcal{R}^\infty_{p,\omega}(0)\, \Longleftrightarrow \, \omega \in B_p(0)\cap B_p(\infty)$.
\end{theo}

\noindent
{\bf Proof.} For each of the statements we only have to prove the implication ``$\Longrightarrow$''.

\medskip
\noindent
``$\Longrightarrow :$''

\noindent {\bf i){:}} As $D_0^{1,p}(\rp,\omega) = \mathcal{R}_{p,\omega}^{0}(0)$, therefore $\mathcal{R}^{0}_{p,\omega}(0)$ is closed subset in $D^{1,p}(\rp,\omega)$.  Hence, by statement i) in Theorem \ref{domknietosczbiorowzero}, we get
$\omega \in B_p(0)$. We have only two possibilities:  a)
$ \omega \in B_p(0)\setminus B_p(\infty)$ or  b) $ \omega \in B_p(0)\cap B_p(\infty)$.
 We will show that condition b) cannot hold.

\noindent  We argue by contradiction.  If the  condition b) was true then, by  Theorem \ref{glowne}, it would imply  $D_0^{1,p}(\rp,\omega) = \mathcal{R}_{p,\omega}^{0}(0)\cap \mathcal{R}_{p,\omega}^{\infty}(0)${,} and consequently
 $\mathcal{R}_{p,\omega}^{0}(0)= \mathcal{R}_{p,\omega}^{0}(0)\cap \mathcal{R}_{p,\omega}^{\infty}(0)$. It would follow that $\mathcal{R}_{p,\omega}^{0}(0)\setminus \mathcal{R}_{p,\omega}^{\infty}(0)=\emptyset$, while the function
\begin{eqnarray}\label{func1}
u(x):= \left\{
\begin{array}{ccc}
0 & {\rm if} & x\in (0,1)\\
x -1  & {\rm if} & x\in [1,2)\\
1 & {\rm if} & x\ge 2
\end{array}
\right.
\end{eqnarray}
belongs to $\mathcal{R}_{p,\omega}^{0}(0)\setminus \mathcal{R}_{p,\omega}^{\infty}(0)$.
We arrive at contradiction, therefore only the condition a) can be true. This proves the statement i).

\medskip
\noindent {\bf ii){:}} We argue similarly as before. As we have $D^{1,p}_0(\rp,\omega)= \mathcal{R}^{\infty}_{p,\omega}(0)$, therefore $\mathcal{R}^{\infty}_{p,\omega}(0)$ is closed.
  Hence, by statement ii) in Theorem \ref{domknietosczbiorowzero}, we get
$\omega \in B_p(\infty)$.

We have only two possibilities: either a)
$ \omega \in B_p(\infty)\setminus B_p(0)$ or  b) $ \omega \in B_p(0)\cap B_p(\infty)$.
We will show that condition b) cannot hold.

\noindent Indeed, if b) would hold, then, from Theorem \ref{glowne}, statement iv), we would deduce that
 $\mathcal{R}^{\infty}_{p,\omega}(0)= \mathcal{R}^{0}_{p,\omega}(0)\cap \mathcal{R}^{\infty}_{p,\omega}(0)$. This would  imply  $\mathcal{R}_{p,\omega}^{\infty}(0)\setminus \mathcal{R}_{p,\omega}^{0}(0)=\emptyset$, while the function

\begin{eqnarray}\label{func2}
u(x):= \left\{
\begin{array}{ccc}
1 & {\rm if} & x\in (0,1)\\
-x +2  & {\rm if} & x\in [1,2)\\
0 & {\rm if} & x\ge 2
\end{array}
\right.
\end{eqnarray}
belongs to $\mathcal{R}_{p,\omega}^{\infty}(0)\setminus \mathcal{R}_{p,\omega}^{0}(0)$.
The contradiction shows that the condition a) holds and this completes the proof of the statement ii).

\medskip
\noindent {\bf iii){:}} By contradiction{,}  let us assume that
$D_0^{1,p}(\rp, \omega)= D^{1,p}(\rp, \omega)$ and
 $\omega \in B_p(0) \cup B_p(\infty)$. Then we have either: a) or b) or c), where  a) $ \omega \in B_p(0)\setminus B_p(\infty)$,  b)\, $ \omega \in B_p(\infty)\setminus B_p(0)$,   c)\, $ \omega \in B_p(0)\cap B_p(\infty)$.

\noindent According to Theorem \ref{glowne}, we then get either $D^{1,p}_0(\rp,\omega)=\mathcal{R}_{p,\omega}^{0}(0)$,\, or $ D^{1,p}_0(\rp,\omega)=\mathcal{R}_{p,\omega}^{\infty}(0)$, or $ D^{1,p}_0(\rp,\omega)=\mathcal{R}_{p,\omega}^{0}(0)\cap \mathcal{R}_{p,\omega}^{\infty}(0)$, respectively.
{Consequently,} we would have either $ D^{1,p}(\rp,\omega)=\mathcal{R}_{p,\omega}^{0}(0)$, or $ D^{1,p}(\rp,\omega)=\mathcal{R}_{p,\omega}^{\infty}(0)$, or
$ D^{1,p}(\rp,\omega)=\mathcal{R}_{p,\omega}^{0}(0)\cap \mathcal{R}_{p,\omega}^{\infty}(0)$, respectively. However, those identities cannot be true. For example, the function $u\equiv 1$ belongs to
$ D^{1,p}(\rp,\omega)$, while it does not belong to any of the sets: $\mathcal{R}_{p,\omega}^{0}(0)$, $\mathcal{R}_{p,\omega}^{\infty}(0)$, $\mathcal{R}_{p,\omega}^{0}(0) \cap \mathcal{R}_{p,\omega}^{\infty}(0)$. \,
 The contradiction proves iii).

\medskip
\noindent {\bf iv){:}} Let us suppose that the implication does not hold, that is  $ D^{1,p}_0(\rp,\omega)= \mathcal{R}^0_{p,\omega}(0)\cap \mathcal{R}^\infty_{p,\omega}(0)$, but
$\omega \not\in B_p(0) \cap B_p(\infty)$.
Then either a) or b) or c) holds, where  a) $ \omega \in B_p(0)\setminus B_p(\infty)$, \, b) $ \omega \in B_p(\infty)\setminus B_p(0)$, \, c) $ \omega \not\in B_p(0)\cup B_p(\infty)$.

\noindent By Theorem \ref{glowne},  those conditions imply: $ \mathcal{R}_{p,\omega}^{0}(0)\cap \mathcal{R}_{p,\omega}^{\infty}(0)=\mathcal{R}_{p,\omega}^{0}(0)$,\, or $ \mathcal{R}_{p,\omega}^{0}(0)\cap \mathcal{R}_{p,\omega}^{\infty}(0)=\mathcal{R}_{p,\omega}^{\infty}(0)$, or $ \mathcal{R}_{p,\omega}^{0}(0)\cap \mathcal{R}_{p,\omega}^{\infty}(0)=D^{1,p}(\rp,\omega)$, respectively.

\noindent In first two situations we would have $\mathcal{R}_{p,\omega}^{0}(0)\setminus \mathcal{R}_{p,\omega}^{\infty}(0)=\emptyset$ or $\mathcal{R}_{p,\omega}^{\infty}(0)\setminus \mathcal{R}_{p,\omega}^{0}(0)=\emptyset$,  which are false, thanks to either (\ref{func1}) or (\ref{func2}), respectively.  The third situation cannot be true also, because the function $u\equiv 1$ belongs to $ D^{1,p}(\rp,\omega)$, but it does not belong to $\mathcal{R}_{p,\omega}^{0}(0)\cap \mathcal{R}_{p,\omega}^{\infty}(0)$. In any case we get the contradiction, which proves  the validity of iv) and ends the proof of the statement.
\hfill$\Box$

\section{Applications}\label{appli}

Let us present several example applications of our results.

\subsection{Application to Hardy inequality}\label{hardyin}

Let us consider classical Hardy and {the} conjugate Hardy operators, respectively:
\begin{equation}\label{hardypresentation}
Hv(t):=\int_0^t v(t)\, dt,\ v\in L^1_{loc}([0,\infty)),\ \ \ \
H^*v(t):= \int_t^\infty v(t)\, dt,\ v\in L^1_{loc}((0,\infty ]).
\end{equation}

\noindent
In Theorems \ref{r-sets} and \ref{r-setsinfty}
 we have shown that Hardy type operators are isometric embeddings between $L^p(\rp ,\omega)$ and $\mathcal{R}^{0}_{p,\omega}(0)$ or $\mathcal{R}^{\infty}_{p,\omega}(0)$, respectively.  Precisely, it follows from Theorem \ref{r-sets} and \ref{r-setsinfty}   that
\begin{eqnarray*}
&H :& L^p(\rp, \omega )\stackrel{{\rm isometry}}{\rightarrow} (\mathcal{R}^{0}_{p,\omega}(0), \|\cdot\|^{(0)}_{{D^{1,p}(\rp,\omega)}}),  \ {\rm   when}  \ \omega\in B_p(0),\\
&H^*:& L^p(\rp, \omega )\stackrel{{\rm isometry}}{\rightarrow} (\mathcal{R}^{\infty}_{p,\omega}(0), \|\cdot\|_{{D^{1,p}(\rp,\omega)}}^{(\infty)}), \ {\rm when}\   \omega\in B_p(\infty).
\end{eqnarray*}
In the first case the inverse is $u\mapsto u^{'}$, while in second case it is $-u^{'}$.

Such identification can be further used to obtain the extended variants of Hardy type inequality, where the class of admissible functions is defined in terms of limits of $u$ at $0$ or at $\infty$.

The necessary and sufficient conditions for boundedness of Hardy operator $H$ and conjugate Hardy operator $H^*$ as acting from $ L^p (\rp, \omega )$ to
$L^q(\rp, h)${, where $1<p,q<\infty$} are known, see e.g. \cite{kuf-opic-book}, \cite{ma}, \cite{muc},  \cite{sin}. For readers convenience we enclose them in the Appendix in Theorems \ref{hardybound} and \ref{hardybound1}.  Let us call them $(C)$ - in case of conditions for $H$, and $(C^*)$ - in case of conditions for $H^*$, respectively.

 We have the following example statement, which deepens
 {our understanding of the Hardy inequality.} As the $B_p(0)$ condition {seems not known before,} in our opinion the result is new.

\begin{theo}[analysis of Hardy inequality]
Suppose that the pair of weight functions $(h,\omega)$, with positive $\omega\in C(\rp)$, satisfies the condition $(C)$ is in Thorem \ref{hardybound},  $1<q, p<\infty$.

Then the following statemets hold.
 \begin{description}
 \item[i)] We have
$h\in L^1_{loc}((0,\infty])$ and $\omega\in B_p(0)$. In particular:
\begin{itemize}
\item the operator
{$Tr^{0}(u):=\displaystyle{\lim_{t\to 0} u(t)}$} is well defined {for} every
{$u\in D^{1,p}(\rp , \omega)$}
\item  the set
 $\mathcal{R}_{p,\omega}^{0}(0)=\{ u\in D^{1,p}(\rp , \omega): Tr^{0}(u)=0\}$ is closed subspace in $D^{1,p}(\rp, \omega)$,
 \item {the inequality}
\begin{eqnarray}\label{sinnamon-related}
\left( \int_{\rp} { \lvert  u(t) \rvert }^q h(t)\, dt\right)^{\frac{1}{q}} \le C \left( \int_{\rp} { \lvert u^{'}(t) \rvert }^p \omega (t)\, dt\right)^{\frac{1}{p}},
\end{eqnarray}
{holds} for {every} $u\in \mathcal{R}_{p,\omega}^{0}(0)$, with {constant} $C {>0}$ independent on $u$.
\end{itemize}
 \item[ii)] When $h\not\in L^1(\rp)$, then inequality \eqref{sinnamon-related}
   with right hand side finite holds precisely on the set $\mathcal{R}^{0}_{p,\omega}(0)$. In particular the embedding
 $X\subseteq L^q (\rp,h)$ cannot be extended to any larger subspace $X$ of $D^{1,p}(\rp,\omega)$.
 \item[iii)] When $h\in L^1(\rp)$, then inequality (\ref{sinnamon-related})
   extends to inequality
 \begin{eqnarray}\label{sinnamon-related-extended}
\| u\|_{L^q(\rp,h)} \le C\left( \| u^{'}\|_{L^p(\rp,\omega)} +  \lvert  Tr^{0} (u) \rvert \right), \  u\in D^{1,p}(\rp , \omega),
\end{eqnarray}
with constant $C$ independent on $u$. Moreover, (\ref{sinnamon-related})
holds on every subspace $V\subseteq D^{1,p}(\rp , \omega)$ which does not contain the nonzero constant functions. In particular
(\ref{sinnamon-related}) holds on  $\mathcal{R}_{p,\omega}^{0}(0)$,   but
 $\mathcal{R}_{p,\omega}^{0}(0)$ is not maximal subspace of $D^{1,p}(\rp , \omega)$ admitted for the validity of (\ref{sinnamon-related}).
 \end{description}
 \end{theo}

\smallskip
\noindent
{\bf Proof.}
{{\bf i):}}
The fact that $h\in L^1_{loc}((0,\infty])$ and $\omega\in B_p(0)$ follows from { Conditions (C), see Theorem \ref{sin}}. The fact that
$\mathcal{R}^{0}_{p,\omega}(0)$ is closed subspace in $D^{1,p}(\rp, \omega)$ and that the operator $Tr^{0}{(\cdot)}$ is well defined on $D^{1,p}(\rp, \omega)$ follows from Theorem \ref{r-sets}. By the same theorem,
$\mathcal{R}^{0}_{p,\omega}(0)$  is precisely the range of Hardy transforms under the action of Hardy operator $H$ applied to
$L^p(\rp, \omega)$. Therefore, under the conditions (C),  (\ref{sinnamon-related}) holds on $\mathcal{R}^{0}_{p,\omega}(0)$. Now we prove the remaining statements.

\smallskip
\noindent
{{\bf ii):}}
We already know that (\ref{sinnamon-related}) holds on $\mathcal{R}^{0}_{p,\omega}(0)$. On the other hand, when $u\not\in \mathcal{R}^{0}_{p,\omega}(0)$ then $Tr^{0}(u) =c\neq 0$, and in such case left hand side in (\ref{sinnamon-related}) cannot be finite as $h$ is not integrable near zero.

\smallskip
\noindent
{{\bf iii):}}
By triangle inequality
$\| u\|_{L^q (\rp ,h)}\le \| u-Tr^{0}(u)\|_{L^q (\rp ,h)} +  \lvert  Tr^{0}(u) \lvert_{L^q (\rp ,h)}$ and {by} the application of (\ref{sinnamon-related}) to $u-Tr^{0}(u)\in \mathcal{R}^{0}_{p,\omega}(0)$, we  easily get (\ref{sinnamon-related-extended}).
Let $V$ be any subspace of $D^{1,p}(\rp,\omega)$ which does not contain
nonzero constant functions. Then the seminorm $\| u\|^{*}_{{D^{1,p}(\rp,\omega)}}$ is the norm
on $V$ and so $X_1:= (V,\|\cdot\|^*_{{D^{1,p}(\rp,\omega)}})$, as well as $X_2:= (V,\|\cdot\|^{(0)}_{{D^{1,p}(\rp,\omega)}})$ (see Theorem \ref{r-sets}),  are Banach spaces. Moreover,  the identity operator $id: X_2\rightarrow X_1$ is continuous linear bijection.
Let us apply
Banach Inverse Mapping Theorem (\cite{rudin}):


\begin{theo}[Banach's Inverse Mapping Theorem].
Let $X,Y$ be Banach spaces and let $T:X\mapsto Y$ be a linear bounded operator. If $T$ is bijective, then $T^{-1}: Y\mapsto X$ is bounded.
\end{theo}

It guarantees that the inverse, $id: X_1\rightarrow X_2$ is  bounded. This implies that the norms $\|\cdot\|^*_{{D^{1,p}(\rp,\omega)}}$ and
$\|\cdot\|^{(0)}_{{D^{1,p}(\rp,\omega)}}$ are comparable on $V$. Thus we can substitute the norm
$\|\cdot\|^{(0)}_{{D^{1,p}(\rp,\omega)}}$ by the norm $\|\cdot\|^*_{{D^{1,p}(\rp,\omega)}}$ in (\ref{sinnamon-related-extended}), when dealing with $u\in V$, with the eventual change of constant in the estimate. Therefore second statement in ii) follows. Last statement follows because, for example, we can consider
$V:= \{ u\in D^{1,p}(\rp,\omega): Tr^{0}(u) =0 \ {\rm or}\ u(1)=0\}$, which is proper subspace of $D^{1,p}(\rp,\omega)$ and a proper
superspace of $\mathcal{R}^{0}_{p,\omega}(0)$.\hfill$\Box$

\begin{rem}[possible analysis of {conjugate} Hardy inequality]\rm~
Similar considerations based on the analysis of Hardy conjugate transform $H^*$ lead to the validity of  (\ref{sinnamon-related}),  equipped with the condition $\displaystyle{\lim_{t\to \infty} u(t)=0}$,  under {Conditions $(C^*)$} for the admitted weights (see Theorem \ref{dualhardy}).
\end{rem}

\subsection{Application to formulation of Dirichlet boundary conditions for solutions of ODE's}\label{bvpb}

The following remark  contributes to the interpretation and well-posedness of boundary conditions of Dirichlet type, in various problems dealing with ODE's.

\begin{rem}\rm
Our analysis allows to interpret precisely Dirichlet type boundary conditions for
$u\in D^{1,p}(\rp,\omega)$  with $1<p<\infty$:
\begin{eqnarray}\label{limits}
\lim_{t\to 0} u(t) &=& c\  \  \hbox{\rm when} \ \omega\in B_p(0),\ \ \ {\rm or}\ \ \
\lim_{t\to \infty} u(t) = c\  \  \hbox{\rm when} \ \omega\in B_p(\infty).
\end{eqnarray}
We already know (see Theorems: \ref{r-sets} and \ref{r-setsinfty}) that, for  $u\in D^{1,p}(\rp,\omega)$, in both cases the above conditions can be equivalently stated as
\begin{eqnarray}\label{boudaryprecise}
\lim_{t\to 0} \frac{u(t)-c}{\left( \int_0^t \omega (\tau )^{-\frac{1}{p-1}  }d\tau\right)^{1-\frac{1}{p}}} &=& 0\  \hbox{\rm when} \ \omega\in B_p(0),\\
\lim_{t\to \infty} \frac{u(t)-c}{\left( \int_t^\infty \omega (\tau )^{-\frac{1}{p-1}  }d\tau\right)^{1-\frac{1}{p}}} &=&0\  \hbox{\rm when} \ \omega\in B_p(\infty).\nonumber
\end{eqnarray}
As in later conditions  the denominators converge to zero, (\ref{boudaryprecise}) is stronger than (\ref{limits}) if we do not assume that
$u\in D^{1,p}(\rp,\omega)$. We can now confirm that the boundary conditions defined by (\ref{boudaryprecise}) are well posed and equivalent for functions in the respective Dirichlet space $D^{1,p}(\rp,\omega)$.
\end{rem}

\subsection{Generalization  of Morrey's inequality}\label{morrey-general}
Morrey's inequality in $1$-dimension says that when  $1<p<\infty$ then
for any $u\in D^{1,p}(\rp {,\omega\equiv 1})$
\begin{equation}\label{morrey}
\| u\|_{C^{0,1-\frac{1}{p}}(\rp)} ={\rm sup}_{x,y\in\rp}\frac{\lvert  u(x)-u(y) \rvert }{ { \lvert x-y \rvert}^{1-\frac{1}{p}}} \le \| u^{'}\|_{L^p(\rp)},
\end{equation}
see e.g \cite{adams-fournier}, Lemma 4.28 on page 99.

\begin{rem}\rm
Using simple modification of inequalities (\ref{a}), we deduce
that when $\omega\in C(\rp)$ {$\omega>0$},  $1<p<\infty$, then
 for any $u\in D^{1,p}(\rp ,\omega)$
\begin{equation*}\label{morrey-general-1}
\| u\|_{{
C_{\omega}^{0,1-\frac{1}{p}}(\rp)
}} :={\rm sup}_{x,y\in\rp}
\frac{ \lvert  u(x)-u(y) \rvert }{d_{w,p}(x,y)^{1-\frac{1}{p}}} \le \| u^{'}\|_{L^p(\rp, \omega)},
\end{equation*}
 where
 $
d_{w,p}(x,y):= \int_x^y \omega^{-\frac{1}{p-1}} (\tau)\ d\tau $ replaces ${\rm dist}(x,y)^{1-\frac{1}{p}}= \lvert  x-y \rvert $.
Observe that the function $d_{\omega,p}(x,y)$ obeys the properties of distance function {on $\rp$} {and replaces $\lvert x-y \rvert$ in \eqref{morrey}}.
\end{rem}

\bigskip

\subsection{Application to complex interpolation theory for {weighted  Dirichlet spaces}}\label{interpo}
In the paper \cite{cwi-ein}, on page 2434, in third question, the authors have asked about complex interpolation results for
 the weighted homogeneous Sobolev spaces, which {in our setting} we call Dirichlet spaces:
$$
D^{1,p}(U , \omega):=\left\{ u\in L^1_{loc} (U) : \frac{\partial u}{\partial x_i} \in L^p(U, \omega),\ {\rm for}\ i=1,\dots ,n\right\} ,
$$
where $U\subseteq \mathbf{R}^n$ is an open set and $\omega: U\rightarrow \rp$ is a given weight. The authors have focused on the case of $p\in [1,\infty )$, $n=1$, and $U=\mathbf{R}$ for the special class of weights, which satisfy the compact boundedness condition as in Definition 1.3 on page 2383. That condition is satisfied by every positive continuous function defined on $\mathbf{R}$.
In that case the mapping:
\begin{equation}\label{izo}
\psi \mapsto  \int_0^x \psi (t) dt  \ \hbox{\rm and its inverse}\ \phi\mapsto \phi^{'}
\end{equation}
give the isomorphic identification between the two Banach couples
\begin{equation*}\label{couples}
(( D^{1,p_0}(\mathbf{R} , \omega_0), ( D^{1,p_1}(\mathbf{R} , \omega_1))\ {\rm and}\ (( L^{p_0}(\mathbf{R} , \omega_0), ( L^{p_1}(\mathbf{R} , \omega_1)).
\end{equation*}
Let $(X,Y)_\theta$ denote the complex interpolation pair between Banach spaces $X,Y$.
It is deduced from Calder\'on type generalization of Stein-Weiss Theorem, as in \cite{cwi-ein}, in Remark 3.2 on page 2397,  that one has:
$$
((L^{p_0} (\mathbf{R},\omega_0), L^{p_1} (\mathbf{R},\omega_1))_\theta = L^{p_\theta}(\mathbf{R},\omega_\theta), 
$$
 where $\omega_\theta^{\frac{1}{p_\theta}} = \omega_0^{ \frac{1-\theta}{p_0} }\omega_1^{\frac{\theta}{p_1}},$
   $\frac{1}{p_\theta} = \frac{1-\theta}{p_0} +\frac{\theta}{p_1}$.

\noindent
From there it follows that
$$
((D^{1,p_0} (\mathbf{R},\omega_0), D^{1,p_1} (\mathbf{R},\omega_1))_\theta = D^{1,p_\theta}(\mathbf{R},\omega_\theta).$$
The precise arguments are submitted in Section A.4 on pages 2439 and 2440 in \cite{cwi-ein}.

\smallskip
\noindent
Consider now the case of $U=\rp$, $p_0,p_1\in (1,\infty )$ and weights
$\omega_0,\omega_1$  such that either a) or b) holds when
a) $\omega_0,\omega_1\in B_p(0)$ or b) $\omega_0,\omega_1\in B_p(\infty)$
are  positive and continuous, $1<p<\infty$.
In case a), the {mapping} (\ref{izo}), while in case b), the {mapping}
\begin{equation*}\label{izoinf}
\psi \mapsto - \int_x^\infty \psi (t) dt  \, \ \hbox{\rm and its inverse}\ \phi\mapsto \phi^{'}
\end{equation*}
give the isomorphic identification between the two Banach couples:
\begin{equation*}\label{couples1}
(( D^{1,p_0}(\rp , \omega_0), ( D^{1,p_1}(\rp , \omega_1))\ {\rm and}\ (( L^{p_0}(\rp , \omega_0), ( L^{p_1}(\rp , \omega_1)).
\end{equation*}
From there, by the same arguments as in \cite{cwi-ein}, we deduce that
$$
((D^{1,p_0} (\rp,\omega_0), D^{1,p_1} (\rp,\omega_1))_\theta = D^{1,p_\theta}(\rp,\omega_\theta )$$
where  $\omega_\theta^{\frac{1}{p_\theta}} = \omega_0^{ \frac{1-\theta}{p_0} }\omega_1^{\frac{\theta}{p_1}}$, 
$\frac{1}{p_\theta} = \frac{1-\theta}{p_0} +\frac{\theta}{p_1} .$
As local $B_p$ conditions $B_p(0)$ and $B_p(\infty)$ have not been analyzed eariler (see Remark \ref{kuopbp}), in our opinion the result is new.

\section{Perspectives for further development, remarks, and open questions}\label{furtherdev}

Let us collect some remarks, focusing on the link with literature and
further possible extensions.

\begin{rem}[local $B_p$ conditions in literature]\label{kuopbp}\rm
In general the localized $B_p$ conditions {are} missing in the literature. However,
in  \cite{kuf-opic}, having positive almost everywhere weight $\omega$ defined on open set $\Omega\subseteq \mathbf{R}^n$, the authors consider the so-called
``exceptional set''  $M_p(\Omega):= $
$$\left\{ x\in \Omega : \int_{\Omega\cap V(x)} \omega^{-1/(p-1)}(y)dy=\infty, \ \hbox{{\rm for every neighborhood  }}\ V(x)\ {\rm of}\ x\right\} .
$$
In our case $\Omega =\rp$ and $\omega \in B_p(\rp)$, so $M_p(\Omega) =\emptyset$, but the extension of the above definition also to   $x\in \bar{\Omega}$ in place of $x\in\Omega$, would lead in our situation to the validation of conditions $B_p(0)$ and $B_p(\infty)$.
\end{rem}

Our results can be extended further in several directions. Let us propose some of them.

\begin{rem}[possible extensions]\rm~\\
(a) {\it The choice of another domain.}
All the results
that we will stated deal with functions defined on the half line.
However, without major changes in the proofs, one can consider instead any
interval  $(a,b)$ in place of $\rp$. We have
focused on functions defined on $\rp$ to make our presentation simpler.
\\
\noindent
(b) {\it Possible discontinuities inside the interval.} We have assumed in all our statements, that weight function $\omega$ is positive and continuous inside the interval $\rp$. It would be interesting to know how much this assumption can be weakened.
\\
\noindent
(c) {\it Higher order Dirichlet spaces.} Instead of $D^{1,p}(\rp,\omega)$, one could consider for example higher order Dirichlet spaces, for example:
\begin{eqnarray*}
D^{k,p}(\rp, \omega)=\{  u:\rp\rightarrow \mathbf{R}:  u\  \hbox{\rm is locally absolutely continuous on}\ \rp \ {\rm and} \\ \| u^{(k)}\|_{L^p(\rp, \omega)}<\infty \}.\nonumber
\end{eqnarray*}
where $k\in \mathbf{N}$ and $u^{(k)}$ is the distributional derivative of $u$, and ask similar questions.
\\
\noindent
{(d) {\it Fractional order Dirichlet spaces.}
 Instead of $D^{1,p}(\rp,\omega)$, one could consider fractional order Dirichlet spaces, where the  derivative $u^{'}$ is replaced by the fractional one, $u^{(\alpha)}$, where $0<\alpha<1$. For example one can use the Caputo, Riemann-Liouville, or Gr\"unwald-Leitnikov derivatives, as discussed for example in the book \cite{sikorski-book}.
 }
\end{rem}


\begin{rem}[{ Muckenhoupt weights}]\rm
In many papers the authors deal with Muckenhoupt weights.
Let $\omega: \mathbf{R}\rightarrow \rp$ be the Muckenhoupt weight,  $1<p<\infty$.
Then, by definition, $\omega$ satisfies the $A_p$ Muckenhoupt condition (see \cite{muc}):
$$
{\rm sup}_I\left( \int_I \frac{1}{\lvert I \rvert}\int_I \omega (t) dt\right)
\left( \int_I \frac{1}{\lvert I \rvert}\int_I \omega (t)^{-1/(p-1)} dt\right)^{p-1} <\infty ,
$$
where $I$ are intervals in $\mathbf{R}$. However our weights are continuous inside $\rp$, they do not need to satisfy the $A_p$ condition. For example, such a one is $\omega(t)=t^{p-1}$, because $\omega^{-1/(p-1)}$ is not integrable near zero.
\end{rem}

According to the discussion made in Section \ref{interpo}, we address the example open question, which naturally arises from our discussion.

\begin{oq}\rm
Suppose that $\omega_0, \omega_1$ are continuous   positive weights defined on
$\rp$,  $1<p<\infty$, such that $\omega_0 \in B_p(0)\setminus B_p(\infty)$ and
$\omega_1 \in B_p(\infty)\setminus B_p(0)$. We ask what is the
complex  interpolation space:
$$
((D^{1,p_0} (\rp,\omega_0), D^{1,p_1} (\rp,\omega_1))_\theta
$$
 where $\theta \in (0,1)$?
\end{oq}

\begin{rem}[similar questions in the Sobolev space setting]\rm
{ Our results can be linked with recent result by Kaczmarek and second author\footnote{
https://arxiv.org/abs/2204.11583}, where the authors deal with power weighted Sobolev spaces (with $\omega (x) = x^\alpha, \alpha \in \mathbf{R}$)
$$
W^{1,p}(\rp,x^\alpha ):= \{ u\in W^{1,1}_{loc}(\rp) :\| u\|_{W^{1,p}(\rp,x^\alpha )}:= \| u\|_{L^{p}(\rp,x^\alpha )} + \| u^{'}\|_{L^{p}(\rp,x^\alpha )} \} ,
$$
and derive similar results such as e.g. the analysis of trace operator, asymptotic behaviour near endpoints: $0$ and $\infty$, density results, applications to complex interpolation theory. As Sobolev spaces and Dirichlet spaces are not the same,  despite similarities,  our approach requires different analysis and it cannot be considered as direct generalization of results by Kaczmarek and second author.
}
\end{rem}

\section{Appendix}\label{appen}

\subsection{The complementary proofs}

%
%
%
%
%
%
%

{\bf Proof of Fact \ref{norm-dirichlet}.}
The proof  is based on modification of arguments from \cite{kuf-opic}, where Sobolev spaces instead of Dirichlet spaces were considered.

Let $U_n:= \{ u_n+c\}_{c\in\mathbf{R}}$ be the Cauchy sequence in $\tilde{D}^{1,p}(\rp,\omega)$. Then for any fixed $a\in \rp$ the function
$$
v_n:= \int_a^t u_n^{'}(\tau)d\tau
$$
is the representative of each $U_n$ in its class in $\tilde{D}^{1,p}(\rp,\omega)$. As $\{ u_n^{'}\}_{n\in \mathbf{N}}$ is the Cauchy sequence in
$L^p (\rp,\omega)$ - the complete space, so $u_n^{'}\stackrel{n\to\infty}{\rightarrow} g$ in $L^p (\rp,\omega)$ for some
$g\in L^p (\rp,\omega)$.
The argumets as in (\ref{bp-rachunek}) allow to conclude that, in the case of $\omega\in B_p$, we have
 $u_n^{'}\to g$ in $L^1_{loc} (\rp)$, as $n\to\infty$, which gives $v_n\to v:= \int_a^t g(\tau)d\tau \in D^{1,p}(\rp,\omega)$ uniformly on compact sets, when $n\to\infty$. This gives
$$
U_n \stackrel{n\to\infty}{\rightarrow} U:=\{ v+c\}_{c\in \mathbf{R}}\
\hbox{\rm}\ {\rm in}\   \tilde{D}^{1,p}(\rp,\omega),$$
because $\| U_n-U\|^{*}_{{\tilde{D}^{1,p}(\rp,\omega)}}= $
$$
=\| \{ v_n-v+c\}_{c\in\mathbf{R}}\|_{{\tilde{D}^{1,p}(\rp,\omega)}}^*=\|(v_n-v)^{'}\|_{L^p(\rp, \omega)}= \| u_n^{'}- g\|_{L^p(\rp, \omega)}\stackrel{n\to\infty}{\rightarrow} 0.
$$
We have  shown that the space $\tilde{D}^{1,p}(\rp,\omega)$ is complete.\hfill$\Box$

\bigskip
\noindent

\noindent
{\bf Proof of Lemma \ref{lemat1}.}
An easy verification shows that $\phi_{(k,K,a)}$ belongs to the admissible class  for the functional, that is {the non-weighted Sobolev space} $W^{1,p}((k,K))$.
We will show that $\phi_{(k,K,a)}$ is the unique minimizer of (\ref{energy}).
We give the proof for $a=1$, because
$\phi_{(k,K,a)}(t)=a\phi_{k,K,1}(t)$. As $E(\cdot)$ is convex functional and the admissible subset of $W^{1,p}((k,K))$ is convex closed set, Direct Methods in the Calculus of Variations (see e.g. \cite{dacorogna}),   give existence of unique minimizer
{of} $E(\cdot)$. Let us call such a minimizer $\phi_{0}$.
Let $T(t):= \displaystyle{\frac{t-K}{k-K}}$ and
\begin{eqnarray*}
\tilde{E}(v):= \int_{(k,K)}{ \lvert \left( T(t) + v(t)\right)^{'} \rvert }^p \omega (t) dt = \int_{(k,K)}{ \lvert  \frac{1}{k-K} + v^{'}(t) \rvert }^p \omega (t) dt,\\
\ {\rm where} \ \ v\in W^{1,p}((k,K),\omega), \ v(k)=v(K)=0.
\end{eqnarray*}
We have
$$
\phi_0 \ \hbox{\rm is minimizer of}\ E(\cdot) \Longleftrightarrow
v_0:= \phi_0-T  \ \hbox{\rm is minimizer of}\ \tilde{E}(\cdot).
$$
By Direct Methods in Calculus of Variations, because
 the functional is nontrivial, coercive and convex on $W^{1,p}((k,K))$ and $\omega\sim 1$ on $[k,K]$, we deduce that there exists a unique
 minimizer of $\tilde{E}$.
To find the minimizer, we compute Euler-Lagrange equation corresponding to the minimizer.

For any $v\in C_0^\infty ((k,K))$ we have
\begin{eqnarray*}
0 &=&{\frac{d}{ds}\tilde{E} (v_0+sv) {\rvert}_{s=0}} \\
&=&
p\int_{(k,K)} \left\{ {\lvert  \frac{1}{k-K} + v_0^{'}(t) \rvert }^{p-1}{\rm sgn}\left( \frac{1}{k-K} + v_0^{'}(t)\right) \omega(t)\right\} v^{'}(t)dt.
\end{eqnarray*}
{As $\omega$ is continuous}, the function inside brackets $\{\cdot\}$ is integrable over $(k,K)$ and its
 weak derivative is zero. Thus this function is constant and hence

\begin{eqnarray}\label{pierwsza-pomoc}
{\lvert  \frac{1}{k-K} + v_0^{'}(t) \lvert }^{p-1}\left\{ ({\rm sgn}\left( \frac{1}{k-K} + v_0^{'}(t)\right) \right\} = \frac{Const}{\omega(t)}.
\end{eqnarray}
Denote $\Phi_r(a):= { \lvert a \rvert}^{r-1} a= { \lvert a \rvert}^r{\rm sign}\, a$, where $r>0, a\in \mathbf{R}$. Then $\Phi_r$ is invertable and
$\Phi_r^{-1}(a)=\Phi_{\frac{1}{r}}(a)$. Applying $\Phi_{\frac{1}{p-1}}$ to both sides in  \eqref{pierwsza-pomoc}, we get
$$
\frac{1}{k-K} + v_0^{'}(t) =\Phi_{\frac{1}{p-1}}\left( \frac{Const}{\omega(t)}\right) = \frac{c_1}{\omega(t)^{\frac{1}{p-1}}},\ \ \hbox{for some}\  c_1\in\mathbf{R}.
$$
This implies
\begin{eqnarray*}
   v_0^{'}(t) &=& -  \frac{1}{k-K}+ c_1\omega(t)^{-1/(p-1)},\\
v_0(t) &=& \int_k^t  v_0^{'}(t) dt =  -  \frac{1}{k-K}(t-k) + c_1 \int_k^t \omega(t)^{-1/(p-1)}dt.
\end{eqnarray*}
for some constant $c_1$. Recalling that $v_0(K)=0$, we {deduce that}
$$c_1 =-\left( \int_{(k,K)} \omega(t)^{-1/(p-1)}dt \right)^{-1},$$
which allows to conclude the statement.
\hfill$\Box$

\subsection{Some results about Hardy and Hardy conjugate operators}

Recall that we deal with Hardy and {conjugate} Hardy operators as in \eqref{hardypresentation}.

\begin{theo}[Conditions (C), \cite{sin}]\label{sin}
For positive weights\\ {$\omega, h:\rp\rightarrow [0,\infty) \cup \{ \infty\}$ and $1<p,q<\infty$,} the  Hardy operator
\begin{equation}\label{hardybound}
H: L^p (\rp, \omega )\rightarrow
L^q(\rp, h)
\end{equation}
{is bounded  if and only if} i) or ii) holds where
\begin{description}
\item[i)] $1<p\le q<\infty$ and
\begin{eqnarray*}
E_1:= {\rm sup}_{t\in (0,\infty)} \left( \int_t^\infty h(s) ds  \right)^{\frac{1}{q}} \left(
\int_0^t \omega (s)^{-\frac{1}{p-1}}ds\right)^{\frac{1}{p^{'}}} <\infty ,
\end{eqnarray*}
\item[ii)]  $1<q< p<\infty$ and
\begin{eqnarray*}
E_2&:=& \int_0^\infty \left( \int_t^\infty h(s) ds  \right)^{\frac{p}{p-q}} \left(
\int_0^t \omega (s)^{-\frac{1}{p-1}}ds\right)^{\frac{p(q-1)}{p-q}} \omega (t)^{-\frac{1}{p-1}}dt  <\infty ,\\
E_3&:=&  \int_0^\infty \left( \int_t^\infty h(s) ds  \right)^{\frac{q}{p-q}} h (t) \left(
\int_0^t \omega (s)^{-\frac{1}{p-1}}ds\right)^{\frac{q(p-1)}{p}} dt <\infty .
\end{eqnarray*}
\end{description}
\end{theo}

\begin{theo}[Conditions $(C^*)$, Theorems: 7.4 and 7.6 in \cite{kuf-opic-book}]\label{dualhardy}
For positive weights $\omega, h: \rp \rightarrow [0,\infty) \cup\{ +\infty\}$ and $1<p,q<\infty$, the conjugate Hardy operator
\begin{equation}\label{hardybound1}
H^*: L^p (\rp, \omega )\rightarrow
L^q(\rp, h)
\end{equation}
is bounded  if and only if
\begin{description}
\item[i)] $1<p\le q<\infty$  and
\begin{eqnarray*}
A&:=& {\rm sup}_{t\in\rp} A(t)<\infty ,  \  {\rm and }\ \lim_{t\to 0} A(t) =\lim_{t\to \infty} A(t) =0,\ {\rm  where}\\
A(t)&:=& \left( \int_0^t h(s) ds  \right)^{\frac{1}{q}} \left(
\int_t^\infty \omega (s)^{-\frac{1}{p-1}}ds\right)^{1-\frac{1}{p}}
.
\end{eqnarray*}
\item[ii)]  $1<q< p<\infty$ and
\begin{eqnarray*}
A&:=& \int_0^\infty \left( \int_0^t h(s) ds  \right)^{\frac{p}{p-q}} \left(
\int_t^\infty \omega (s)^{-\frac{1}{p-1}}ds\right)^{\frac{p(q-1)}{p-q}} \omega (t)^{-\frac{1}{p-1}}dt .
\end{eqnarray*}
\end{description}

\end{theo}

\section*{Acknowledgements} 
 The work of CC has been partially supported by Istituto Nazionale di Alta Matematica/Gruppo Nazionale per l'Analisi Matematica, la Probabilit\'a e le loro Applicazioni. A.K. wish to thank 
  several very pleasant stays at Istituto per le Applicazioni del Calcolo "Mauro Picone"
Consiglio Nazionale delle Ricerche, where our colaboration originated.
  
\section{Declarations}

\subsection*{Ethical Approval}
Not aplicable.
 
\subsection*{Competing interests} 
The authors have no competing interests as defined by Springer, or other interests that might be perceived to influence the results and/or discussion reported in this paper.

\subsection*{Authors' contributions}
The authors contributed equally to this work.

\subsection*{Fundings} 
The work of CC has been partially supported by Istituto Nazionale di Alta Matematica/Gruppo Nazionale per l'Analisi Matematica, la Probabilit\'a e le loro Applicazioni. 
 
\subsection*{Availability of data and materials} 
All of the material in the manuscript is owned by the authors and no permissions are required.


\bigskip

\begin{thebibliography}{99}

\bibitem{adams-fournier}  ADAMS R.,  FOURNIER J. F.: {\rm
Sobolev spaces,} Second edition. Pure and Applied Mathematics (Amsterdam), 140. Elsevier/Academic Press, Amsterdam (2003).

\bibitem{cwi-ein}  CWIKEL M.,  EINAV A.: {\em
Interpolation of weighted Sobolev spaces,}
J. Funct. Anal. {\bf  277}(7), 2381--2441 (2019).

\bibitem{dacorogna}  DACOROGNA B.:
{\em Direct methods in the calculus of variations,}
Second edition. Applied Mathematical Sciences, 78. Springer, New York (2008).

\bibitem{haj-kalam}  HAJ\L{}ASZ P.,  KA\L{}AMAJSKA A.: {\em
Polynomial asymptotics and approximation of Sobolev functions,}
Studia Math. {\bf 113}(1), 55--64 (1995).



\bibitem{akkppstudia09}  KA\L{}AMAJSKA A.,  PIETRUSKA-PA\L{}UBA K.:
{\em On a variant of the Hardy inequality between weighted Orlicz spaces,}
Studia Math. {\bf 193}(1), 1--28 (2009).

\bibitem{kjf}  KUFNER A.,  JOHN O.,  FUCIK S.: {\em Function spaces},  Monographs and Textbooks on Mechanics of Solids and Fluids, Mechanics: Analysis. Noordhoff International Publishing, Leyden; Academia, Prague (1977).

\bibitem{kuf-opic-book}  KUFNER A.,  OPIC B.:
{\em Hardy-type inequalities,}
Pitman Research Notes in Mathematics Series, 219. Longman Scientific $\&$ Technical, Harlow (1990).

\bibitem{kuf-opic}  KUFNER A.,  OPIC B.: {\em How to define reasonably weighted sobolev spaces,} Commentationes
Mathematicae Universitatis Carolinae, {\bf 25}(3), 537--554 (1984).


\bibitem{ma}  MAZ'YA V.G.: {\em Sobolev Spaces}, Springer, Berlin (1985).

\bibitem{sikorski-book}  MEERSCHAERT M.M.,  SIKORSKII A.:
{\em Stochastic models for fractional calculus,} Second edition. De Gruyter Studies in Mathematics, 43. De Gruyter, Berlin (2019).

\bibitem{muc}  MUCKENHOUPT B.:  {\em Weighted norm inequalities for the Hardy maximal function,} Transactions of the American Mathematical Society {\bf  165}, 207--226 (1972).

\bibitem{rudin}  RUDIN W.:
{\em Functional analysis,}
Second edition. International Series in Pure and Applied Mathematics. McGraw-Hill, Inc., New York (1991).



\bibitem{sin}  SINNAMON G.J.:
{\em Weighted Hardy and Opial-type inequalities,}
J. Math. Anal. Appl. {\bf 160}(2),  434--445  (1991).



\bibitem{sobolev}  SOBOLEV S.L.: {\em
The density of compactly supported functions in the space $L_p^{(m)}(E)$,} ( in Russian), Sibirsk. Mat. \v{Z}. {\bf 4}, 673--682 (1963).



\end{thebibliography}
\end{document}